\documentclass[reqno,12pt]{amsart}
\usepackage{amscd,amsfonts,mathrsfs,amsthm,amssymb,amsmath,mathtools}
\usepackage{upgreek}
\usepackage{csquotes}
\usepackage{xcolor}
\usepackage{enumerate}
\usepackage{bbm}
\usepackage{multicol}
\usepackage{color}
\usepackage{array}
\usepackage{ae}
\usepackage{accents}
\usepackage{epsfig}
\usepackage[e]{esvect}
\usepackage[nobysame]{amsrefs}
\usepackage[margin=1in]{geometry}
\usepackage[toc,title,page]{appendix}
\usepackage{hyperref}
\hypersetup{colorlinks=true,linkcolor=blue,filecolor=blue,urlcolor=blue, citecolor=blue}

\usepackage{ragged2e}
\usepackage{marginfix}
\let\oldtocsection=\tocsection
\let\oldtocsubsection=\tocsubsection
\let\oldtocsubsubsection=\tocsubsubsection
\renewcommand{\tocsection}[2]{\hspace{0em}\oldtocsection{#1}{#2}}
\renewcommand{\tocsubsection}[2]{\hspace{1em}\oldtocsubsection{#1}{#2}}
\renewcommand{\tocsubsubsection}[2]{\hspace{2em}\oldtocsubsubsection{#1}{#2}}

\def\image #1#2#3#4{
	\begin{figure}[!ht]
		\centering
		\begin{minipage}[t]{0.02\textwidth}
		\end{minipage}
		\begin{minipage}[t]{#4\textwidth}       
			\includegraphics[width=\textwidth]{#1.eps}
		\end{minipage}%
		\caption{#2}
		\label{#3}
	\end{figure}
}

\def\real     #1{{\mathbb R^{#1}}}
\def\complex  #1{{\mathbb C^{#1}}}
\def\natural  #1{{\mathbb N^{#1}}}

\def\lap      {\Delta }
\def\dt       {\partial_{t}}

\def\equationcolor {\color{black}}
\def\textcolor     {\color{black}}

\def\bcoleq    {\begin{equation}\equationcolor}
\def\ecoleq    {\textcolor\end{equation}}
\def\bcoleqn   {\equationcolor\begin{eqnarray}}
\def\ecoleqn   {\end{eqnarray}\textcolor}

\def\y{\mathrm{u}}

\def\det{\operatorname{det}}
\def\R{\mathbb{R}}

\newcommand{\la}{\langle}
\newcommand{\ra}{\rangle}
\newcommand{\na}{\nabla}

\newcommand{\De}{\Delta}

\let\epsilon\varepsilon
\let\rho\varrho

\newtheorem{theorem}{Theorem}[section]
\newtheorem{mythm}{Theorem}

\newtheorem*{question*}{Question}

\newtheorem{lemma}[theorem]{Lemma}

\newtheorem{corollary}[theorem]{Corollary}

\theoremstyle{definition}
\newtheorem*{assumption*}{$\lambda_{1}$-Condition}

\newtheorem{remark}[theorem]{Remark}

\usepackage{tikz}
\usetikzlibrary{patterns}
\newcommand{\hatchedsquare}{%
  \begin{tikzpicture}[scale=0.25, baseline=-0.3ex]
    \draw[pattern=north east lines] (0,0) rectangle (1,1);
  \end{tikzpicture}%
}

\newenvironment{proofof}[1]{%
  \par\noindent\textbf{Proof of #1.}\enspace\ignorespaces
}{%
  \hfill $\hatchedsquare$\par\medskip
}

\usepackage{titlesec}
\titleformat{\section}
  {\normalfont\large\bfseries\centering} 
  {\thesection.}                       
  {.3em}                                 
  {}

\usepackage{amsthm}
\usepackage{amssymb} 

\numberwithin{equation}{section}

\begin{document}

\title[Ancient solutions of MCF in a slab]{Ancient solutions of mean curvature flow in a slab}
\author[N.M. M{\o}ller ]{\textsc{Niels Martin M{\o}ller}}
\author[A. Savas-Halilaj]{\textsc{Andreas Savas-Halilaj}}
\author[K. Smoczyk]{\textsc{Knut Smoczyk}}

\address{Niels Martin M{\o}ller\newline
{Department of Mathematical Sciences,
Copenhagen Centre for Geometry and Topology\\
University of Copenhagen, DK-2100 Copenhagen, Denmark} \newline
{\sl E-mail:} {\bf nmoller@math.ku.dk}
}

\address{Andreas Savas-Halilaj\newline
{Department of Mathematics,
Section of Algebra \!\&\! Geometry, \!\!
University of Ioannina,
45110 Ioannina, Greece} \newline
{\sl E-mail:} {\bf ansavas@uoi.gr}
}

\address{Knut Smoczyk\newline
Institut f\"ur Differentialgeometrie and
Riemann Center for Geometry \& Physics,
Leibniz Universit\"at Hannover, Welfengarten 1,
30167 Hannover, Germany\newline
{\sl E-mail address:} {\bf smoczyk@math.uni-hannover.de}
}

\renewcommand{\subjclassname}{  \textup{2020} Mathematics Subject Classification}
\thanks{N.M. M\o{}ller was supported by DFF Sapere Aude 7027-00110B, by Carlsberg Foundation Semper Ardens CF21-0680 and by CPH-GEOTOP-DNRF151. K. Smoczyk was funded by the Deutsche Forschungsgemeinschaft (DFG, German Research Foundation) -- 569348513 (SM 78/8-1).}
\parindent = 0 mm
\hfuzz     = 6 pt
\parskip   = 3 mm

\begin{abstract}
In this paper, we study complete ancient solutions of the mean curvature flow in arbitrary codimension contained in slabs in Euclidean space. Under a natural principal-angle condition, we prove that a translating soliton contained in a slab of width less than $\pi$ cannot be bounded from below in the translating direction, and characterize the grim reaper cylinder in form of a rigidity result at the critical width $\pi$. We derive corresponding results for ancient solutions confined in slabs.
\end{abstract}

\setcounter{tocdepth}{1}

\maketitle

\section{Introduction}

An ancient solution $F:M^m\times(-\infty,T)\to\real{m+k}$ ($T\le\infty$) of the mean curvature flow is said to be {\em confined in a slab} if its image is contained within the region bounded by two parallel hyperplanes. Since the mean curvature flow is invariant under ambient isometries, we can assume without loss of generality throughout this work that the slab is given by $[-a/2,a/2]\times\real{m+k-1}$ for some $a>0$, which we refer to as the {\em width} of the slab. 

Translating solitons constitute one of the basic classes of self-similar solutions to the mean curvature flow, arising naturally as models for Type-II singularities. Specifically, an immersion $F:M^m\to\real{m+k}$
is a {\em translating soliton} if its mean curvature vector satisfies
$$
\mathbf H=V^\perp,
$$
for some constant non-zero vector $V\in\real{m+k}$. Throughout the paper, we
normalize the translation vector by $|V|=1$. Moreover, we will always assume our submanifolds to be connected. Since every translating soliton generates an eternal (and hence ancient) solution to the mean curvature flow, translators serve as a fundamental source of examples and models for ancient solutions.

A particularly interesting class of translating solitons and ancient solutions consists of those confined in 
slabs. The simplest example is the {\em grim reaper curve}, which is contained in a slab of width $\pi$, together 
with its products with Euclidean factors. Over the years, considerable attention has been devoted to the 
existence, classification, and rigidity of such solutions, particularly regarding translating solitons in 
codimension one.
A first rigidity result in this direction was obtained by
Mart\'in, P\'erez-Garc\'ia, Savas-Halilaj, and Smoczyk \cite{MSS1},
who proved that a connected properly embedded translator in
$\mathbb{R}^3$ with uniformly bounded genus on compact sets,
which is $C^1$-asymptotic outside a cylinder to two planes,
must be either a plane or the grim reaper cylinder.
This result was subsequently extended to arbitrary dimensions
by Gama and Mart\'in in \cite{GM}, who characterized vertical
hyperplanes and tilted grim reaper cylinders among translators
that are $C^1$-asymptotic outside a non-vertical cylinder
to two half-hyperplanes.

The existence and classification of {\em graphical} translators
in slabs have also been investigated extensively.
Hoffman, Ilmanen, Mart\'in, and White \cite{HIMW} obtained
a classification of complete graphical translators in
$\mathbb{R}^3$, including existence and uniqueness results
for strictly convex translating graphs over strips of
width greater than $\pi$. In arbitrary dimensions, Bourni, Langford and Tinaglia
\cite{BLT} constructed strictly convex translators contained in slabs of every width greater than $\pi$, and showed that no such translator can exist in a slab of width less than or equal to $\pi$. On the other hand, Chini and
M{\o}ller \cite{CM} established bi-halfspace and convex-hull
theorems for properly immersed translators, providing further
insight into the possible global configurations of these
submanifolds and, in particular, the occurrence of slabs.
More recently, Gama, Mart\'in, and M{\o}ller \cite{GMM22}
investigated translators of finite genus and finite entropy
contained in slabs in $\mathbb{R}^3$, obtaining an entropy
formula in terms of the number of asymptotic wings,
as well as a uniqueness result for translating pitchforks.
Subsequently, in \cite{GMM25}, they developed a structure theory for this class of translators, showing that their ends are asymptotic to unique collections of planes and possibly tilted grim reaper cylinders, alongside further rigidity results for translators bounded from below.

In a related direction, Impera, M{\o}ller, and Rimoldi
\cite{impera} established rigidity and non-existence
results for finite-entropy translators contained in slabs in $\R^3$,
including a characterization of the grim reaper cylinder
at the critical width $\pi$.
We note that in higher dimensions ($m \ge 4$), such examples exist trivially via 
cylindrical extensions. If $N^{m-1} \subset \real{m}$ is a minimal hypersurface, the 
product $M^m = N^{m-1} \times \real{} \subset \real{m+1}$ forms a codimension-one 
translator in the direction of the factor $\real{}$. Jorge and Xavier \cite{JX} 
constructed examples of complete nonplanar and non-proper minimal surfaces contained 
between two parallel planes in $\real{3}$. While Hoffman and Meeks excludes properly immersed nonplanar minimal surfaces contained in a slab in $\real{3}$, $n$-dimensional catenoids $C^{n}\subset\real{n+1}$, for $n \ge 3$, provide proper minimal embedded hypersurfaces in a slab. Rescaling these higher-dimensional catenoids generates codimension-1 translators of dimension $m=n+1 \ge 4$ contained in slabs of arbitrarily small width.
Motivated by these developments, the following question arises naturally (which has now been answered in the negative for $\mathbb{R}^3$ in \cite{GMM26}, assuming finite topology and finite entropy):

\begin{question*}
Under which conditions can one classify  non-planar proper translators 
in a slab, depending on the width and the codimension?
\end{question*}

A first aim of the present paper is therefore to investigate translating solitons contained in slabs in arbitrary codimension. In this setting, the geometry is substantially more subtle than in the hypersurface case. Indeed, if $V$ denotes the unit translating direction and $W$ is a constant unit vector orthogonal to $V$, then the fixed two-plane
$$
P\doteq\operatorname{span}\{V,W\}
$$
interacts with both the tangent and the normal spaces of the submanifold. In codimension one, the normal space is one-dimensional, rendering this interaction largely degenerate. In higher codimension, however, the normal projection of the plane $P$ may be genuinely two-dimensional.

As an instructive example to illustrate these complications, let us consider a naive attempt to extend the results in \cite{CM} to higher codimension. A natural analog of the wedge theorem would assert that no properly immersed translator $\Sigma^m\subset\mathbb{R}^{m+k}$ for $k\ge2$ lies in the intersection of two halfspaces with linearly independent normals $W_1,W_2\perp V$. 

However, this fails, as shown by e.g. $\sqrt{2}\Gamma\times\sqrt{2}\Gamma\subset\mathbb{R}^4$ with coordinates $(x_1, y_1, x_2, y_2)$, where $\Gamma=\{(x,-\log\cos x):|x|<\pi/2\}$. This is a codimension-two translator with velocity $V=(E_{y_1}+E_{y_2})/\sqrt{2}$, which nevertheless lies strictly inside the vertical wedge
$x_1, x_2>-\pi/\sqrt{2}$. To see why attempts to extend the argument fail for this example, recall that the 
strategy in the  hypersurface proof is to apply the Omori-Yau maximum principle to a distance function, yielding a contradiction 
with its strict subharmonicity at infinity. In the higher-codimension example, both factors become asymptotically 
vertical along the maximizing sequence, making both horizontal directions asymptotically normal to the surface. 
The otherwise positive Laplacian term then tends to zero, preventing any contradiction from arising.
In contrast, a hypersurface's one-dimensional normal space cannot accommodate both directions, forcing the 
Laplacian of the distance function to remain positively lower-bounded and thereby ensuring that the non-existence 
result holds. One therefore sees that angle conditions are expected to enter in higher codimension, to replace 
the role played by the ambient wedge geometry in codimension-one; see also Lemma \ref{prop:P-intersection}.

This naturally leads to the investigation of the Gauss map $\mathcal G:M\to {\rm Gr}_m(\real{m+k})$, given by
$$
\mathcal G(p)\doteq T_pM,
$$
and to the relative position of the tangent space $T_pM$ and the fixed $2$-plane
$P=\operatorname{span}\{V,W\}.$
Since $P$ is 2-dimensional, this relative position is described by two
principal angles
$$
0\leq\beta_1\leq\beta_2\leq{\pi}/{2}.
$$
These angles measure how the plane $P$ is distributed between the
tangent and normal spaces of the translator and provide the natural
geometric condition needed in our arguments. We refer to
Section \ref{Gauss} for the precise definition and the corresponding
projection formulas.

Having identified the relevant geometry of the Gauss map, we introduce
the two natural height functions
$$
\y\doteq\langle V,F\rangle
\quad\text{and}\quad
\alpha\doteq\langle W,F\rangle.
$$
The function $\y$ measures the height in the translating direction,
while $\alpha$ measures the position of the translator in the direction
orthogonal to the slab. In particular, containment in a slab is
equivalent to boundedness of $\alpha$, and the width of the smallest
such slab is $\operatorname{osc}(\alpha)$.
In our main results, we shall assume that the principal angles satisfy
\begin{equation}\label{b1b2intro}
\beta_1+\beta_2\leq{\pi}/{2},\tag{$\beta$}
\end{equation}
and that the slab width, measured by the function $\alpha$, satisfies
\begin{equation}\label{anglea}
\operatorname{osc}(\alpha)<\pi.
\end{equation}
These two assumptions arise naturally in a number of geometric settings.
As a matter of fact, the principal-angle condition \eqref{b1b2intro} is automatic for
hypersurfaces and for coisotropic submanifolds of $\complex{n}$, hence
also for Lagrangian translators. In the graphical setting, it is
equivalent to the weak area-decreasing condition. Thus the geometric
hypothesis appearing in our main theorem is satisfied by several natural classes of translating solitons.

The first main result of the paper is the following:

\begin{mythm}\label{thma}
Let $F:M^m\to\real{m+k}$, $m\geq2$, be a properly immersed translating soliton moving in the unit direction $V$. Let
$W\in\real{m+k}$ be a constant unit vector orthogonal to $V$, and set
$\y\doteq\langle V,F\rangle$ and $\alpha\doteq\langle W,F\rangle.$
Assume that the principal angles between the plane $P\doteq{\rm span}\{V,W\}$ and the tangent spaces
of the translator satisfy the angle condition \eqref{b1b2intro}.
Suppose further that $\operatorname{osc}(\alpha)<\pi.$
Then the height function $\y$ cannot be bounded from below.
\end{mythm}

In this borderline situation, we prove a rigidity theorem showing that,
under a natural extremal assumption, the translator splits as the product
of a grim reaper curve with a complete minimal submanifold.

\begin{mythm}\label{thmb}
Let $F:M^m\to\real{m+k}$, $m\geq2$, be a properly immersed translating soliton
moving in the unit direction $V$. Let
$W\in\real{m+k}$ be a constant unit vector orthogonal to $V$, and set
$\y\doteq\langle V,F\rangle$ and $\alpha\doteq\langle W,F\rangle.$
Assume that the principal angles between the plane $P\doteq{\rm span}\{V,W\}$ and the tangent spaces
of the translator satisfy the angle condition \eqref{b1b2intro}.
Suppose further that $\operatorname{osc}(\alpha)=\pi$. Then the following statements are equivalent:
\begin{enumerate}[\normalfont(i)]
    \item The non-negative function $e^u \cos(\alpha-\alpha_0)$, $\alpha_0=\sup\alpha-\pi/2$, attains a local minimum.
    \item $e^u \cos(\alpha-\alpha_0)$ is a positive constant.
    \item $M$ is a Riemannian product $\Gamma \times \Sigma$, where $\Gamma \subset P$ is a $1$-dimensional grim reaper curve and $\Sigma \subset \mathbb{R}^{m+k-2}$ is a complete $(m-1)$-dimensional minimal submanifold.
\end{enumerate}
\end{mythm}

Slab confinement also plays an important role in the study of
ancient solutions to mean curvature flow. In this direction,
Bourni, Langford and Tinaglia \cite{BLTancient} constructed
compact convex ancient solutions, known as ancient pancakes,
which lie in slabs of critical width $\pi$.
Chini and M{\o}ller \cites{CMancient,CM} established bi-halfspace
and convex-hull theorems for ancient mean curvature flows,
providing further restrictions on their global geometry, while Ma \cite{MaHalfspace} constructed examples of ancient solutions that sweep out a half-space.
In a related development, Bourni, Langford and Tinaglia
\cite{BLTpoly} investigated ancient mean curvature flows
arising from polytopes, constructing new examples and studying
their asymptotic geometry in relation to grim hyperplanes.
More recently, Bourni, Langford and Lynch \cite{BLL}
characterized collapsing convex ancient solutions in terms
of their asymptotic behavior, showing that collapsing is
equivalent to being asymptotic to at least one grim hyperplane.
Motivated by these developments, we also investigate ancient
solutions confined in slabs in arbitrary codimension. We obtain the result that ancient solutions confined in a slab determined by the direction $W$ entirely sweep out in time all coordinate directions $V$ perpendicular to $W$ for which our angle condition \eqref{b1b2intro} holds - in particular, all coordinate directions $V$ perpendicular to $W$ in the case of hypersurfaces.

\begin{mythm}\label{thmc}
Let $F:M^m\times(-\infty,T)\to\mathbb R^{m+k}$ be a complete ancient solution to mean curvature flow. Assume that, for each
$[a,b]\subset(-\infty,T)$, the map \((x,t)\mapsto(F(x,t),t)\) is proper on 
$M\times[a,b]$.
Let $V, W \in \mathbb{R}^{m+k}$ be constant orthonormal vectors, 
and set $u \coloneqq \langle V, F \rangle$ and $\alpha \coloneqq \langle W, F \rangle$. 
Assume that the principal angles between $P \coloneqq \operatorname{span}\{V, W\}$ and the tangent spaces 
of the solution satisfy the condition \eqref{b1b2intro} 
on $M \times (-\infty, T)$. If the solution is confined to a slab, i.e., if $\operatorname{osc}(\alpha) \le L$ for some constant $L > 0$, 
then $u$ can be neither uniformly bounded from above nor from below on $M \times (-\infty, T)$. In particular, in the codimension one case ($k=1$), any coordinate function $u$ associated with a vector $V$ orthogonal to $W$ automatically satisfies the angle condition \eqref{b1b2intro}, and thus must possess this property.
\end{mythm}

While we assume properness in our theorems to ensure the applicability of the Omori-Yau maximum principle, our results hold under any alternative condition that validates the principle. Sufficient geometric and analytic conditions for this principle, including curvature bounds, volume growth conditions, and the existence of suitable exhaustion functions, can be found in the work of Pigola, Rigoli and Setti \cite{PRS}. For ancient solutions, the corresponding parabolic Omori-Yau maximum principle is required; for more details see Section \ref{OY}.

\noindent{}\textbf{AI Disclosure.}
During the preparation of this work, the authors used AI tools solely for editorial assistance, 
including refining the presentation and language, and converting text into \LaTeX{} code. All 
mathematical calculations, proofs, and arguments were independently derived and checked by the 
authors, who assume full responsibility for the content of the article.
 
 \section{Examples of translators confined in slabs}
 
In this section, we provide examples of translating solitons that are 
contained in slabs. We begin with the grim reaper cylinders and their
tilted versions, and then recall further codimension-one examples. We also discuss the examples of
Joyce-Lee-Tsui Lagrangian translators, which 
exhibit a rather different phenomenon, since they can be contained in slabs of 
arbitrarily small width.

\subsection*{Grim reaper cylinder.}
The simplest example of a translating soliton contained in a slab is the
{\em grim reaper cylinder}; see Figure \ref{fig: GrimReaper}. After normalizing the translating vector to have unit
length, the standard grim reaper curve $\Gamma\subset\real{2}$ is given by
$$
\Gamma\doteq\left\{(\alpha,-\log\cos \alpha):-{\pi}/{2}<\alpha<{\pi}/{2}\right\},
$$
and hence $\Gamma\times\mathbb{R}^{m-1}$
is contained in a slab of width $\pi$.
\image{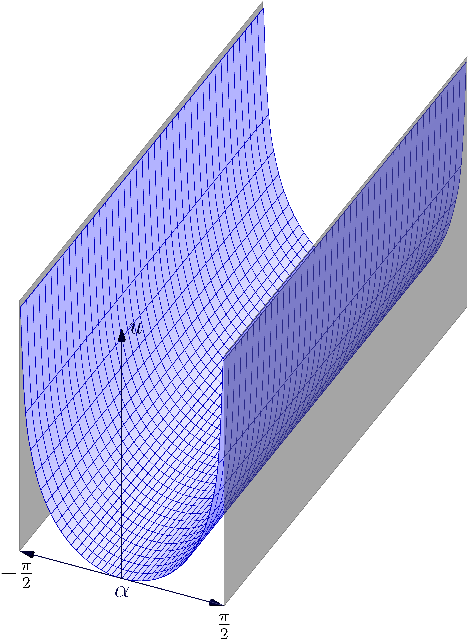}{The $2$-dimensional grim reaper cylinder $\Gamma\times\real{}\subset\real{3}$.}{fig: GrimReaper}{.35}

\subsection*{Tilted grim hyperplanes.}\label{tilted}
Let $V=(0,\dots,0,1)\in\R^{n+1}$ be the unit translating vector in
$\mathbb{R}^{n+1}$. Fix
$a=(a_2,\ldots,a_n)\in\mathbb{R}^{n-1},$
and consider the graphical translator, moving in the direction $V$, described by the map
$$
F(x_1,\ldots,x_n)=(x_1,\ldots,x_n,u(x_1,\ldots,x_n)),
$$
where $u$ is the function
$$
u(x_1,\ldots,x_n)=-(1+|a|^2)
\log\Big(\cos\tfrac{x_1}{\sqrt{1+|a|^2}}\Big)+\sum_{j=2}^n a_jx_j,
$$
defined on the slab
$$
\Big(-\frac{\pi}{2}\sqrt{1+|a|^2},\frac{\pi}{2}\sqrt{1+|a|^2}\Big)\times\mathbb{R}^{n-1}.
$$
For $a=0$ this reduces to the standard grim reaper cylinder
$\Gamma\times\mathbb{R}^{n-1},$ whose slab has width $\pi$.
If $a\neq0$, the slab width is strictly larger than $\pi$. Let us mention that Bourni, Langford, and Tinaglia \cite{BLT} proved that in every dimension $n\geq2$, any strictly convex graphical translator in codimension one in $\real{n+1}$ over $V^\perp$, normalized to move with unit speed $V$, cannot be contained in a slab of width less than or equal to $\pi$; see also \cite{HIMW}.

\subsection*{$\Delta$-wings.}
Another important family is formed by the $\Delta$-{\em wings} constructed by
Hoffman, Ilmanen, Mart\'in, and White \cite{HIMW}. In dimension two, for every $w>\pi$
there exists a complete strictly convex translating graph over a strip of
width $w$. Higher-dimensional analogues over slabs of prescribed width
$w>\pi$ were also constructed.

\subsection*{Pitchfork translators.}
Hoffman, Mart\'in, and White constructed in \cite{HMW1} non-graphical complete properly
embedded translators in $\mathbb{R}^3$, known as {\em pitchfork translators},
which are contained in slabs. These examples are simply connected and have
several planar ends.

\subsection*{Annuloids.}
Hoffman, Mart\'in, and White \cites{HMW3,HMW2} also constructed remarkable families of complete properly
embedded translating annuli, usually called {\em annuloids}, which are contained
in slabs. These examples provide non-simply-connected translating solitons
with several planar ends.

\subsection*{Product examples in higher codimension.}
Examples in arbitrary codimension can be obtained by taking products. More
precisely, if $F_1:M_1\to\mathbb{R}^{n_1}$
is a translating soliton with translating vector $V$ and
$F_2:M_2\to\mathbb{R}^{n_2}$
is a minimal immersion, then
$$
F=(F_1,F_2):M_1\times M_2\to\mathbb{R}^{n_1+n_2}
$$
is again a translating soliton, with translating vector $(V,0)$.
If $F_1(M_1)$ is contained in the slab
$$
S_{a,b}=\left\{x\in\mathbb{R}^{n_1}:a<\langle x,W\rangle<b\right\},
$$
where $W\in\mathbb{R}^{n_1}$ is a fixed unit vector satisfying
$\langle W,V\rangle=0$, then the product translator
$$
F=(F_1,F_2):M_1\times M_2\to\mathbb{R}^{n_1}\times\mathbb{R}^{n_2}
$$
is contained in the slab
$$
S_{a,b}\times\mathbb{R}^{n_2}
=\left\{(x,y)\in\mathbb{R}^{n_1}\times\mathbb{R}^{n_2}:
a<\langle x,W\rangle<b\right\}.
$$
In particular, taking the product of the grim reaper curve with any complete
minimal submanifold produces translating solitons of arbitrary codimension
contained in a slab of width $\pi$. Moreover, Nadirashvili \cite{nadi} constructed complete minimal surfaces contained in a ball in $\mathbb{R}^3$. After rescaling, such a surface may be contained in a ball of arbitrarily small radius. Taking its product with a line and choosing the translating direction along the line factor, we obtain complete translating solitons contained in slabs of arbitrarily small width.

\subsection*{The examples of Joyce-Lee-Tsui.}
We conclude with some very interesting examples of Lagrangian translating solitons
constructed by Joyce-Lee-Tsui in \cite{JLT}*{Corollary I}. It turns out that these examples may be
contained in slabs of arbitrarily small width.
Given $\tau>0$ and
$a_1,\ldots,a_{n-1}>0$, define
$$
P(t)\doteq\frac{1}{t^2}\left(e^{\tau t^2}\prod_{k=1}^{n-1}(1+a_kt^2)-1\right)
$$
and, for $j\in\{1,\ldots,n-1\}$, define
$$
\phi_j(y)\doteq\int_0^y\frac{dt}{\left(\frac{1}{a_j}+t^2\right)\sqrt{P(t)}}.
$$
Joyce-Lee-Tsui prove in \cite{JLT}*{Corollary I} that the map $F:\real{n}\to\complex{n}$
given by
$$
F(x_1,\dots,x_{n-1},y)=\left(z_1(x_1,\dots,x_{n-1},y),\dots,z_n(x_1,\dots,x_{n-1},y)\right),
$$
where
$$
z_j
=
x_j\sqrt{\frac{1}{a_j}+y^2}\,\,\,e^{i\phi_j(y)},
\quad\text{for}\,\,\,
j\in\{1,\ldots,n-1\},
$$
and
$$
z_n
=
\frac{1}{2}y^2
-\frac{1}{2}\sum_{j=1}^{n-1}x_j^2
-\frac{i}{\tau}\sum_{j=1}^{n-1}\phi_j(y)
-\frac{i}{\tau}
\arg\left(y+iP(y)^{-1/2}\right),
$$
parametrizes a closed embedded Lagrangian submanifold of
$\complex{n}$, diffeomorphic to $\real{n}$, which is a translating
soliton with translating vector
$$
V=(0,\ldots,0,\tau)\in\complex{n},
$$
where the last entry denotes the real direction of the last complex
coordinate.

We now take $\tau=1$. Thus the translating vector is the unit vector $V=(0,\ldots,0,1)$.
Writing $z_n=x_n+iy_n$, we have $V=\partial_{{x_n}}$ and $W=JV=\partial_{y_n}.$
In particular,
\begin{equation}\label{joyceheight}
\y\doteq \langle V,F\rangle=\operatorname{Re}z_n
=\frac{1}{2}y^2-\frac{1}{2}\sum_{j=1}^{n-1}x_j^2
\quad\text{and}\quad\alpha\doteq\langle W,F\rangle=\operatorname{Im}z_n.
\end{equation}
If $\vartheta$ denotes the Lagrangian angle, then
Joyce-Lee-Tsui prove in \cite{JLT}*{Proposition 2.5} that every
connected Lagrangian translating soliton with translating vector $V$
satisfies
$$
\vartheta=-\langle JV,F\rangle+c
$$
for some constant $c$. For the above parametrization, one may take
$c=0$, and hence, when $\tau=1$, we have that
\begin{equation}\label{thetaa}
\vartheta=-\langle W,F\rangle=-\alpha.
\end{equation}
More explicitly,
$$
\vartheta(y)=\sum_{j=1}^{n-1}\phi_j(y)+\arg\left(y+iP(y)^{-1/2}\right).
$$
The size of the oscillation is determined in
\cite{JLT}*{Corollary I}. More precisely, there exist numbers
$$
\bar\phi_1,\ldots,\bar\phi_{n-1}\in\left(0,{\pi}/{2}\right)
$$
such that
$$
\lim_{y\to+\infty}\phi_j(y)=\bar\phi_j,\quad\lim_{y\to-\infty}\phi_j(y)=-\bar\phi_j
\quad\text{and}\quad\sum_{j=1}^{n-1}\bar\phi_j<\frac{\pi}{2}.
$$
Moreover, Joyce-Lee-Tsui prove in the same corollary that, for fixed
$\tau>0$, the map
$$
(a_1,\ldots,a_{n-1})
\longmapsto
(\bar\phi_1,\ldots,\bar\phi_{n-1})
$$
is a one-to-one correspondence onto the set
$$
\left\{(\bar\phi_1,\ldots,\bar\phi_{n-1})\in\left(0,\frac{\pi}{2}\right)^{n-1}:
\sum_{j=1}^{n-1}\bar\phi_j<\frac{\pi}{2}\right\}.
$$
They also show in \cite{JLT}*{Corollary I} that the Lagrangian angle
varies between
$$
\sum_{j=1}^{n-1}\bar\phi_j
\quad\text{and}\quad
\pi-\sum_{j=1}^{n-1}\bar\phi_j.
$$
Consequently, due to \eqref{thetaa}, we get that
$$
\operatorname{osc}(\vartheta)=\operatorname{osc}(\alpha)=\pi-2\sum_{j=1}^{n-1}\bar\phi_j.
$$
It follows directly from the one-to-one correspondence above that
$$\bar\phi\doteq\sum_{j=1}^{n-1}\bar\phi_j$$
can be chosen arbitrarily close to
$\pi/2$, while keeping $\tau=1$. Therefore, for every
$\varepsilon>0$, there exists a unit-speed Joyce-Lee-Tsui translator
such that
$$
\operatorname{osc}(\alpha)<\varepsilon.
$$
Equivalently, these examples are contained in slabs orthogonal to
$W=JV$ of arbitrarily small width. Notice that Joyce-Lee-Tsui
explicitly make this observation in the paragraph immediately following
\cite{JLT}*{Corollary I}: using Proposition~2.5, they point out that when
the oscillation of the Lagrangian angle is small compared with the
length of the translating vector, the coordinate
$\operatorname{Im}z_n$ is confined to a small interval. In their
notation, this interval is centered near the hyperplane
$\operatorname{Im}z_n=-{\pi}/(2\tau).$
It is also worth observing that the height function
$$
\y=\frac{1}{2}y^2-\frac{1}{2}\sum_{j=1}^{n-1}x_j^2
$$
given in \eqref{joyceheight} is unbounded from below.

\section{The Gauss map and the principal angles}\label{Gauss}
Let $F:M^m\to\real{m+k}$
be an immersed translating soliton with unit translating direction $V$,
so that
$$
\mathbf H=V^\perp.
$$
The Gauss map $\mathcal G:M\to {\rm Gr}_m(\real{m+k})$ of $F$ is defined to be
$$
\mathcal G(p)\doteq dF(T_pM)\,\,\,(\text{which for simplicity we identify with}\,\,\, T_pM).
$$
Fix a unit vector $W\in\real{m+k}$ orthogonal to $V$, and consider the
two-dimensional plane
$$
P\doteq \operatorname{span}\{V,W\}.
$$
The relative position of $P$ and $dF(T_pM)$ is described by their principal
angles. More precisely, let $\pi_p^\top:P\to T_pM$ be the orthogonal projection and consider the self-adjoint endomorphism $A_p$ given by
$$
A_p\doteq(\pi_p^\top)^*\pi_p^\top:P\to P.
$$
If $1\geq\mu_1(p)\geq\mu_2(p)\geq 0$ are the eigenvalues of $A_p$, we define the angles $\beta_1(p)$
and $\beta_2(p)$ by the relations
$$
0\leq\beta_1(p)\leq\beta_2(p)\leq{\pi}/{2}
\quad\text{and}\quad \mu_i(p)=\cos^2\beta_i(p),
\quad i\in\{1,2\}.
$$
Thus $\cos\beta_1$ and $\cos\beta_2$ are the singular values of
$\pi_p^\top$. When $m=1$, the map $\pi_p^\top$ has rank at most one and
hence $\beta_2=\pi/2$. With this convention, the two angles are defined
in every dimension. With the above ordering,
$\beta_1$ and $\beta_2$ also define continuous functions on $M^m$.
Note that the angles $\beta_1$ and $\beta_2$ depend only on the value
of the Gauss map and measure the relative position of
the fixed plane $P$ and the tangent space of the translator. For this reason, we shall refer to them as the {\em principal angles of the immersion}.

The following elementary identities will be useful.

\begin{lemma}\label{lemma:principal_angles}
At every point $p\in M^m$, we have
$$
|V^\top\wedge W^\top|^2=\cos^2\beta_1\cos^2\beta_2
\quad\text{and}\quad |V^\perp\wedge W^\perp|^2=\sin^2\beta_1\sin^2\beta_2.
$$
Moreover,
$$
|V^\top|^2+|W^\top|^2
=
\cos^2\beta_1+\cos^2\beta_2.
$$
\end{lemma}

\begin{proof}
Since $\{V,W\}$ is an orthonormal basis of $P$, the matrix of $A_p$
with respect to this basis is the Gram matrix of $V^\top$ and
$W^\top$. Hence
$$
\det A_p=|V^\top|^2|W^\top|^2-\langle V^\top,W^\top\rangle^2=|V^\top\wedge W^\top|^2.
$$
Since the eigenvalues of $A_p$ are
$\cos^2\beta_1$ and $\cos^2\beta_2$, it follows that
$$
|V^\top\wedge W^\top|^2=\cos^2\beta_1\cos^2\beta_2.
$$
Similarly, the Gram matrix of $V^\perp$ and $W^\perp$ is $I-A_p$.
Its eigenvalues are $\sin^2\beta_1$ and $\sin^2\beta_2$, and therefore
$$
|V^\perp\wedge W^\perp|^2=\sin^2\beta_1\sin^2\beta_2.
$$
Finally,
$$
|V^\top|^2+|W^\top|^2=\operatorname{tr}A_p=\cos^2\beta_1+\cos^2\beta_2.
$$
This completes the proof.
\end{proof}

The geometric condition appearing in our main results has several
equivalent formulations.

\begin{lemma}\label{lemma:angle_equivalences}
The following conditions are equivalent:
\begin{enumerate}[\rm(1)]
\item $\beta_1+\beta_2\leq{\pi}/{2}.$
\medskip
\item $|V^\perp\wedge W^\perp|\leq |V^\top\wedge W^\top|.$
\medskip
\item $|V^\top|^2+|W^\top|^2\geq1.$
\end{enumerate}
\end{lemma}

\begin{proof}
By Lemma~\ref{lemma:principal_angles}, we have
$$
|V^\top\wedge W^\top|^2-|V^\perp\wedge W^\perp|^2=\cos^2\beta_1+\cos^2\beta_2-1.
$$
Since
$$
\cos^2\beta_1+\cos^2\beta_2-1=\cos(\beta_1+\beta_2)\cos(\beta_1-\beta_2)
\quad\text{and}\quad \cos(\beta_1-\beta_2)\geq 0,
$$
the first two conditions are equivalent. The equivalence with
$$|V^\top|^2+|W^\top|^2\geq 1$$
follows from Lemma~\ref{lemma:principal_angles}. This completes the proof.
\end{proof}

We now record several geometric situations in which the condition
$$
\beta_1+\beta_2\leq{\pi}/{2}
$$
is automatically satisfied.

\subsection*{A tangential direction in $\boldsymbol{P}$.}

The simplest sufficient condition is that the plane $P$ has a
nontrivial intersection with the tangent space.

\begin{lemma}\label{prop:P-intersection}
If $P\cap T_pM\neq\{0\},$ then $\beta_1(p)=0$ and hence
$$
\beta_1(p)+\beta_2(p)\leq {\pi}/{2}.
$$
In particular, the principal-angle assumption appearing in our main theorem imposes
no additional restriction for hypersurfaces.
\end{lemma}
\begin{proof}
If $P\cap T_pM\neq\{0\}$, there exists a unit vector
$X\in P\cap T_pM$. Thus $|\pi_p^\top(X)|=1,$
and hence the largest singular value of $\pi_p^\top$ is equal to one.
Therefore $\cos\beta_1(p)=1$, so that $\beta_1(p)=0$.
For hypersurfaces the preceding condition is always satisfied. Indeed, since
$$
\dim P+\dim T_pM=m+2>m+1,
$$
we have $\dim(P\cap T_pM)\geq1$, and the conclusion follows from the preceding argument.
This completes the proof.
\end{proof}
\subsection*{Graphs over the plane $\boldsymbol{P}$.}

The angle condition also admits a natural interpretation for graphical
surfaces. Let $\Omega\subset P$ and let $f:\Omega\to P^\perp$
be a smooth map. Consider its graph $F:\Omega\subset\R^2\to\R^{2+k}$
given by
$$
F(x)=x+f(x).
$$
Let $0\leq\lambda_1\leq\lambda_2$
be the singular values of $df$. Then the principal angles between $P$
and the tangent space to the graph satisfy
$$
\tan\beta_i=\lambda_i,\quad i\in\{1,2\}.
$$
Consequently, $\beta_1+\beta_2\leq{\pi}/{2}$
if and only if $\lambda_1\lambda_2\leq1.$
Since the Jacobian of $f$ satisfies
$$
|\Lambda^2df|=\lambda_1\lambda_2,
$$
we obtain the following characterization.

\begin{lemma}\label{prop:graph_angles}
Let $F:\Omega\subset\R^2\to P\oplus P^\perp$ be the immersion given by $F(x)=x+f(x)$.
Then
$$
\beta_1+\beta_2\leq{\pi}/{2}
$$
if and only if $f$ is weakly area-decreasing, that is, if and only if
$$
|\Lambda^2df|\leq1.
$$
\end{lemma}
\subsection*{Coisotropic submanifolds.}

There is another important class for which our principal-angle condition
is automatic. Let $F:M^m\to\complex{n}$
be an immersed submanifold, where $\complex{n}$ is equipped with its
standard complex structure $J$ and symplectic form $\omega$. Recall that
$M=F(M^m)$ is {\em coisotropic} if
$$
(T_pM)^\omega\subset T_pM,\quad\text{for every}\,\,\, p\in M,
$$
or, equivalently, if and only if,
$$
J(N_pM)\subset T_pM, \quad\text{for every}\,\,\, p\in M.
$$
\begin{lemma}\label{prop:coisotropic_angles}
Let $F:M^m\to\complex{n}$
be a coisotropic translating soliton with unit translating direction
$V$, and set $W=JV$. Then
$$
\beta_1+\beta_2\leq {\pi}/{2}.
$$
\end{lemma}

\begin{proof}
Write $V=V^\top+V^\perp.$
Since $M$ is coisotropic, it follows that $JV^\perp\in T_pM.$
Moreover,
$$
W^\top=(JV^\top)^\top+JV^\perp.
$$
The two terms on the right-hand side are orthogonal, since
$$
\langle JV^\top,JV^\perp\rangle=\langle V^\top,V^\perp\rangle=0.
$$
Consequently, $|W^\top|^2\geq|V^\perp|^2.$
Hence
$$
|V^\top|^2+|W^\top|^2 \geq |V^\top|^2+|V^\perp|^2=1.
$$
Lemma~\ref{lemma:angle_equivalences} now gives the desired result.
\end{proof}

\subsection*{Lagrangian submanifolds.}
Every Lagrangian submanifold is coisotropic, so the previous lemma applies. In fact, in the Lagrangian case the angle condition is always
saturated.

\begin{lemma}\label{prop:lagrangian_angles}
Let $F:M^n\to\complex{n}$
be a Lagrangian translating soliton with unit translating direction
$V$, and set $W=JV$. Then $\beta_1+\beta_2={\pi}/{2}.$
\end{lemma}

\begin{proof}
Since $M$ is Lagrangian, $J(T_pM)=N_pM$ and $J(N_pM)=T_pM$, for any $p\in M$.
Therefore $(JV^\top)^\top=0$ and $W^\top=JV^\perp.$
It follows that $|W^\top|^2=|V^\perp|^2$ and hence
$$
|V^\top|^2+|W^\top|^2=|V^\top|^2+|V^\perp|^2=1.
$$
By Lemma~\ref{lemma:angle_equivalences}, this is precisely the equality
case $\beta_1+\beta_2={\pi}/{2}.$
This completes the proof.
\end{proof}

\section{Gradient relations and drift identities}
We begin with some elementary relations between the gradients of two important
height functions.

\begin{lemma}\label{lem:grad_relations}
Let $F:M^m\to\real{m+k}$ be an immersed submanifold. Let
$V,W\in\real{m+k}$ be fixed orthonormal vectors, and define the
height functions $\y$ and $\alpha$ by
$$
\y\doteq\langle V,F\rangle
\quad\text{and}\quad
\alpha\doteq\langle W,F\rangle.
$$
Then the following hold true:
\begin{enumerate}[\rm (1)]
\item
The intrinsic gradients of $u$ and $\alpha$ satisfy
$$\nabla\y=V^\top\quad\text{and}\quad \nabla\alpha=W^\top,$$
where with $\{\cdot\}^{\top}$ we denote the orthogonal projection on the
tangent bundle of the submanifold.

\medskip

\item The gradients of $u$ and $\alpha$ satisfy the identities,
$$
1-|\nabla\y|^2=|V^\perp|^2
\quad\text{and}\quad
1-|\nabla\alpha|^2=|W^\perp|^2.
$$
Moreover,
$$
\langle\nabla\y,\nabla\alpha\rangle
=
-\langle V^\perp,W^\perp\rangle
$$
and
$$
(1-|\nabla\y|^2)(1-|\nabla\alpha|^2)-\langle\nabla\y,\nabla\alpha\rangle^2={\rm S}(\y,\alpha),
$$
where
$$
{\rm S}(\y,\alpha)\doteq|V^\perp\wedge W^\perp|^2=|V^\perp|^2|W^\perp|^2-\langle V^\perp,W^\perp\rangle^2\geq{0}.
$$

\item The sum of squares of the gradients of $u$ and $\alpha$ satisfies 
$$
|\nabla\y|^2+|\nabla\alpha|^2=1+{\rm C}(\y,\alpha)-{\rm S}(\y,\alpha),
$$
where
$$
{\rm C}(\y,\alpha)\doteq|V^\top\wedge W^\top|^2=|\nabla\y|^2|\nabla\alpha|^2-\langle\nabla\y,\nabla\alpha\rangle^2\geq 0.
$$
\end{enumerate}
\end{lemma}

\begin{proof}
We decompose the constant vectors $V$ and $W$ into their tangential
and normal components as
$$
V=V^\top+V^\perp
\quad\text{and}\quad
W=W^\top+W^\perp.
$$
(1) For every tangent vector $Z\in T_pM$, we have
$$
\langle\nabla\y,Z\rangle=Z(\y)=Z\langle V,F\rangle=\langle V^\top,Z\rangle.
$$
Therefore, $\nabla\y=V^\top.$
In exactly the same way, it follows that $\nabla\alpha=W^\top.$

\noindent (2) Since the tangential and normal components are orthogonal, taking the
squared norm of the first decomposition and using $|V|=1$ gives
\begin{equation}\label{duv}
1=|V|^2=|\nabla\y|^2+|V^\perp|^2.
\end{equation}
Similarly, since $|W|=1$,
\begin{equation}\label{duw}
1=|W|^2=|\nabla\alpha|^2+|W^\perp|^2.
\end{equation}
Since $V$ and $W$ are orthogonal, we also have
$$
0=\langle V,W\rangle=
\langle\nabla\y+V^\perp,
\nabla\alpha+W^\perp\rangle=\langle\nabla\y,\nabla\alpha\rangle
+\langle V^\perp,W^\perp\rangle.
$$
Therefore,
$$
\langle\nabla\y,\nabla\alpha\rangle
=
-\langle V^\perp,W^\perp\rangle.
$$
Combining the above identities, we obtain
$$
(1-|\nabla\y|^2)(1-|\nabla\alpha|^2)
-
\langle\nabla\y,\nabla\alpha\rangle^2
=|V^\perp|^2|W^\perp|^2-\langle V^\perp,W^\perp\rangle^2={\rm S}(\y,\alpha)\ge 0.
$$
(3) Expanding the above identity and using \eqref{duv} and 
\eqref{duw}, we get
\begin{eqnarray*}
{\rm S}(\y,\alpha)&\!\!\!=\!\!\!&1-|\nabla\y|^2-|\nabla\alpha|^2
+|\nabla\y|^2|\nabla\alpha|^2-\langle\nabla\y,\nabla\alpha\rangle^2\\
&\!\!\!=\!\!\!&
1-|\nabla\y|^2-|\nabla\alpha|^2
+{\rm C}(\y,\alpha),
\end{eqnarray*}
from where it follows that
$$
|\nabla\y|^2+|\nabla\alpha|^2
=1+{\rm C}(\y,\alpha)-{\rm S}(\y,\alpha).
$$
This completes the proof.
\end{proof}

\begin{lemma}\label{prop:laplacians}
Let $F:M^m\to\real{m+k}$ be a translator moving in the direction $V$
under mean curvature flow. Then
$$
\lap\y=1-|\nabla\y|^2\quad\text{and}\quad
\lap\alpha=-\langle\nabla\y,\nabla\alpha\rangle.
$$
\end{lemma}

\begin{proof}
For every constant vector $Z\in\real{m+k}$, by the Beltrami identity
we have
$$
\lap\langle Z,F\rangle=\langle Z,\lap F\rangle=\langle Z,\mathbf H\rangle.
$$
Since $\mathbf H=V^\perp$,
it follows that
$$
\lap\langle Z,F\rangle
=
\langle Z,V^\perp\rangle.
$$
Setting $Z=V$ and taking into account Lemma~\ref{lem:grad_relations}(2),
we obtain
$$
\lap\y=\langle V,V^\perp\rangle=|V^\perp|^2=1-|\nabla\y|^2.
$$
Setting $Z=W$ and keeping in mind again Lemma~\ref{lem:grad_relations}(2), we obtain
$$
\lap\alpha=\langle W,V^\perp\rangle=\langle W^\perp,V^\perp\rangle=-\langle\nabla\y,\nabla\alpha\rangle.
$$
This completes the proof.
\end{proof}

Define now the function $f_\sigma:M^m\to\R$ given by
$$
f_\sigma\doteq e^\y\cos(\sigma\alpha),
$$
where $\sigma$ is a constant such that
$$
\cos(\sigma\alpha)>0.
$$

\begin{lemma}\label{prop:f_sigma}
The function $f_\sigma$ satisfies the identity
\begin{equation}\label{eq:f_sigma_raw}
\lap f_\sigma=\langle\nabla f_\sigma,\nabla\y\rangle
+f_\sigma(1-|\nabla\y|^2-\sigma^2|\nabla\alpha|^2).
\end{equation}
Equivalently, we can rewrite the above identity in the form
\begin{equation}\label{thm:drift_identity}
\lap f_\sigma=\langle\nabla f_\sigma,\nabla\y
+\lambda \mathcal{X}\rangle-\lambda f_\sigma
+(1+\lambda)f_\sigma\big({\rm S}(\y,\alpha)-{\rm C}(\y,\alpha)\big),
\end{equation}
where
$$
\mathcal{X}\doteq \nabla\y+\sigma\tan(\sigma\alpha)\nabla\alpha
\quad\text{and}\quad
\lambda=\frac{\sigma^2-1}{1+\sigma^2\tan^2(\sigma\alpha)}.
$$
\end{lemma}

\begin{proof}
From Lemma \ref{prop:laplacians}, we get
$$
\nabla e^\y=e^\y\nabla\y\quad\text{and}\quad
\lap(e^\y)=e^\y(\lap\y+|\nabla\y|^2)=e^\y.
$$
Moreover,
$$
\nabla\cos(\sigma\alpha)
=
-\sigma\sin(\sigma\alpha)\nabla\alpha,
$$
and hence
\begin{eqnarray*}
\lap\cos(\sigma\alpha)
&\!\!\!=\!\!\!&
-\sigma\sin(\sigma\alpha)\lap\alpha
-\sigma^2\cos(\sigma\alpha)|\nabla\alpha|^2\\
&\!\!\!=\!\!\!&
\sigma\sin(\sigma\alpha)
\langle\nabla\y,\nabla\alpha\rangle
-\sigma^2\cos(\sigma\alpha)|\nabla\alpha|^2.
\end{eqnarray*}
Furthermore,
\begin{eqnarray*}
2\langle\nabla e^\y,\nabla\cos(\sigma\alpha)\rangle
&\!\!\!=\!\!\!&
-2\sigma e^\y\sin(\sigma\alpha)
\langle\nabla\y,\nabla\alpha\rangle.
\end{eqnarray*}
Combining these identities, we obtain
$$
\lap f_\sigma=f_\sigma(1-\sigma^2|\nabla\alpha|^2)-\sigma e^\y\sin(\sigma\alpha)\langle\nabla\y,\nabla\alpha\rangle.
$$
On the other hand,
$$
\nabla f_\sigma=e^\y\cos(\sigma\alpha)\nabla\y
-\sigma e^\y\sin(\sigma\alpha)\nabla\alpha
=f_\sigma\nabla\y-\sigma e^\y\sin(\sigma\alpha)\nabla\alpha.
$$
Consequently,
\begin{eqnarray*}
\langle\nabla f_\sigma,\nabla\y\rangle
&\!\!\!=\!\!\!&
f_\sigma|\nabla\y|^2
-\sigma e^\y\sin(\sigma\alpha)
\langle\nabla\y,\nabla\alpha\rangle.
\end{eqnarray*}
Subtracting this identity from the expression for $\lap f_\sigma$
gives
$$
\lap f_\sigma-\langle\nabla f_\sigma,\nabla\y\rangle
=f_\sigma(1-|\nabla\y|^2-\sigma^2|\nabla\alpha|^2),
$$
which proves the identity \eqref{eq:f_sigma_raw}. Because
$$
1+\lambda=\frac{\sigma^2(1+\tan^2(\sigma\alpha))}
{1+\sigma^2\tan^2(\sigma\alpha)}
$$
we have
\begin{eqnarray*}
1-|\nabla\y|^2-\sigma^2|\nabla\alpha|^2
=
(1+\lambda)(1-|\nabla\y|^2-|\nabla\alpha|^2)
-\lambda+\lambda(|\nabla\y|^2-\sigma^2\tan^2(\sigma\alpha)|\nabla\alpha|^2).
\end{eqnarray*}
By Lemma~\ref{lem:grad_relations}(3), we deduce that
\begin{eqnarray*}
1&\!\!\!-\!\!\!&|\nabla\y|^2-\sigma^2|\nabla\alpha|^2\\
&\!\!\!=\!\!\!&-\lambda-(1+\lambda){\rm C}(\y,\alpha)+(1+\lambda)
{\rm S}(\y,\alpha)
+\lambda(|\nabla\y|^2-\sigma^2\tan^2(\sigma\alpha)|\nabla\alpha|^2).
\end{eqnarray*}
Since
$$f_\sigma=e^\y\cos(\sigma\alpha)>0,$$
we have
$$
\log f_\sigma=\y+\log\cos(\sigma\alpha)
\quad\text{and}\quad
\nabla\log f_\sigma=\nabla\y-\sigma\tan(\sigma\alpha)\nabla\alpha.
$$
Hence
$$
\langle\nabla\log f_\sigma,\mathcal{X}\rangle
=
|\nabla\y|^2
-\sigma^2\tan^2(\sigma\alpha)|\nabla\alpha|^2,
$$
from where it follows that
$$
\lambda(|\nabla\y|^2-\sigma^2\tan^2(\sigma\alpha)|\nabla\alpha|^2)
=\frac{1}{f_\sigma}\langle\nabla f_\sigma,\lambda \mathcal X\rangle.
$$
Therefore,
\begin{eqnarray*}
f_\sigma(1-|\nabla\y|^2-\sigma^2|\nabla\alpha|^2)
=-\lambda f_\sigma
+(1+\lambda)f_\sigma\big({\rm S}(\y,\alpha)-{\rm C}(\y,\alpha)\big)
+\langle\nabla f_\sigma,\lambda \mathcal X\rangle.
\end{eqnarray*}
Substituting this into \eqref{eq:f_sigma_raw}, we obtain
the desired identity.
\end{proof}

\section{The Omori-Yau maximum principle}\label{OY}
In this section, we recall the Omori-Yau maximum principle
and its parabolic counterpart, which will play an important role
in the proofs of our main results. We refer to the monograph
of Pigola, Rigoli and Setti \cite{PRS} for a comprehensive
treatment of the Omori-Yau maximum principle and sufficient
geometric and analytic conditions ensuring its validity.

Let $(M,g)$ be a complete Riemannian manifold. We say that
the {\em Omori-Yau maximum principle} for the Laplacian holds on $M$
if, for every function $u\in C^2(M)$ satisfying $\sup u<\infty,$
there exists a sequence of points $p_j\in M$ such that
\begin{itemize}
\item $\lim u(p_j)=\sup u$,
\medskip
\item $\lim |\nabla u|(p_j)= 0$,
\medskip
\item $\limsup\Delta u(p_j)\leq 0$.
\end{itemize}
Applying this principle to $-u$, one obtains the corresponding
statement for functions bounded from below.

The classical theorem of Omori \cite{Omori} and Yau \cite{Yau} guarantees the validity
of this principle on complete Riemannian manifolds whose Ricci
curvature is bounded from below. More generally, Pigola, Rigoli
and Setti \cite{PRS} established sufficient conditions involving
the existence of suitable exhaustion functions, extending
the applicability of the principle beyond curvature assumptions.
In particular, the following result holds:

\begin{lemma}
Let $F:M^m\to\mathbb{R}^{m+k}$ be a properly immersed
submanifold with bounded mean curvature vector.
Then the Omori-Yau maximum principle for the Laplacian
holds on $M$.
\end{lemma}

In particular, every properly immersed translating soliton
satisfies the Omori-Yau maximum principle. Indeed, if the
translating velocity $V$ is a unit vector, then
$\mathbf H=V^\perp$, and hence $|\mathbf H|\leq 1$.
Moreover, finite entropy implies properness for complete translating
solitons, as proved by Chini \cite{Chini20}*{Lemma 24} (see also \cite{ChiniThesis} and Kim and Pyo \cite{kim}*{Theorem 1.2}). Consequently, every
complete finite-entropy translator satisfies the Omori-Yau maximum
principle.

For ancient solutions to mean curvature flow, we shall need
a parabolic version of the Omori-Yau maximum principle. Such a principle
was established by Ma \cite{Ma} for proper mean curvature flows
in ambient manifolds satisfying suitable curvature assumptions.
In the setting of ancient mean curvature flows, Chini and
M{\o}ller \cite{CMancient}*{Lemma 2} obtained a corresponding
maximum principle for properly immersed hypersurfaces in
Euclidean space. The following version of the parabolic Omori-Yau maximum
principle is adapted to our purposes. It follows from the
penalization argument of Chini and M{\o}ller
\cite{CMancient}*{Lemma 2}, suitably modified to account for
an initial time boundary and to allow arbitrary codimension.

\begin{lemma}[Parabolic Omori-Yau maximum principle]
Let $F:M^m\times[t_0,T)\to\mathbb{R}^{m+k}$
be a smooth solution to mean curvature flow satisfying
$\partial_tF=\mathbf H$. Assume that, for every
$s\in(t_0,T)$, the spacetime immersion
$(x,t)\mapsto(F(x,t),t)$
is proper on $M\times[t_0,s]$. Let $u\in C^{2,1}(M\times[t_0,T))$ be bounded from below
and suppose that
$$
u_*\doteq \inf_{M\times[t_0,T)}u <\inf_M u(\cdot,t_0).
$$
Then there exists a sequence $(x_j,t_j)\in M\times(t_0,T)$ such that
$$
\lim u(x_j,t_j)= u_*,\quad \lim |\nabla u|(x_j,t_j)= 0\quad\text{and}\quad 
\limsup(\partial_t-\Delta)u(x_j,t_j)\leq 0.
$$
Here $\nabla$ and $\Delta$ denote the gradient and
the Laplace-Beltrami operator with respect to
the metric induced by $F_t=F(\cdot,t)$.
\end{lemma}

\begin{proof}
By the definition of $u_*$ and the strict inequality
$$u_*<{\inf}_M u(\cdot,t_0),$$
we may choose a sequence $(y_j,\tau_j)\in M\times(t_0,T)$
such that
$$
u_*\leq u(y_j,\tau_j)\leq u_*+1/j.
$$
Set
$$
\varepsilon_j\doteq \frac{1}{j(1+|F(y_j,\tau_j)|^2)}>0.
$$
Then
$$
\varepsilon_j\leq 1/j \quad\text{and}\quad \varepsilon_j|F(y_j,\tau_j)|^2\le 1/j.
$$
Define the sequence of functions
$$
u_j\doteq u+\varepsilon_j|F|^2.
$$
Choose $s_j\in(\tau_j,T)$. By the properness assumption, the function $u_j$ attains
its minimum on $M\times[t_0,s_j]$ at some point $(x_j,t_j)$.
Indeed, since $u\geq u_*$ and $\varepsilon_j>0$, we have
$$
u_j\geq u_*+\varepsilon_j|F|^2.
$$
Hence, every sublevel set of $u_j$ is a closed subset
of the inverse image, under the proper spacetime immersion
$(x,t)\mapsto(F(x,t),t)$, of a compact cylinder
$\overline B_R(0)\times[t_0,s_j]$.
It follows that the sublevel sets of $u_j$ are compact,
and therefore $u_j$ attains its minimum.
By construction,
$$
u_j(x_j,t_j)\leq u_j(y_j,\tau_j) \leq u_*+2/j.
$$
Since $u_*<{\inf}_M u(\cdot,t_0),$ for sufficiently large $j$ we have
$$
u_*+2/j<{\inf}_M u(\cdot,t_0).
$$
On the other hand, $u_j(x_j,t_j)\leq u_*+2/j,$
whereas
$$u_j(x,t_0)\geq u(x,t_0)\geq{\inf}_M u(\cdot,t_0).$$
Consequently, the minimum point $(x_j,t_j)$ cannot
lie on the initial time slice, and hence $t_j>t_0$.
Moreover,
$$
u_*\leq u(x_j,t_j)\leq u_*+2/j,
$$
whence
$$u(x_j,t_j)\to u_*.$$
Moreover, since
$$u(x_j,t_j)\geq u_*,$$
we have
$$
\varepsilon_j|F(x_j,t_j)|^2=u_j(x_j,t_j)-u(x_j,t_j)\leq u_*+2/j-u(x_j,t_j)\leq 2/j.
$$
At the spatial minimum of the function $u_j$, it follows that
$$\nabla u_j(x_j,t_j)=0.$$
Consequently,
$$
|\nabla u|(x_j,t_j)=2\varepsilon_j|F^\top(x_j,t_j)|\leq 2\varepsilon_j|F(x_j,t_j)|\rightarrow 0.
$$
Furthermore, since $t_j>t_0$ and $(x_j,t_j)$ is a
spacetime minimum of $u_j$ on $M\times[t_0,s_j]$, we have
$$
\nabla u_j(x_j,t_j)=0\quad\text{and}\quad \Delta u_j(x_j,t_j)\geq0.
$$
Moreover,
$$\partial_t u_j(x_j,t_j)\leq 0.$$
Indeed, if $t_j<s_j$, the time derivative vanishes,
whereas if $t_j=s_j$, the same inequality follows
by considering the left time derivative at the minimum.
Consequently,
$$
(\partial_t-\Delta)u_j(x_j,t_j)\leq0.
$$
On the other hand, the mean curvature flow equation yields
$$
(\partial_t-\Delta)|F|^2=-2m.
$$
Recalling that
$$u_j=u+\varepsilon_j|F|^2,$$
we obtain
$$
(\partial_t-\Delta)u(x_j,t_j)=(\partial_t-\Delta)u_j(x_j,t_j)+2m\varepsilon_j 
\leq2m\varepsilon_j.
$$
Since $\varepsilon_j\to0$, it follows that
$$
\limsup(\partial_t-\Delta)u(x_j,t_j)\leq0,
$$
which completes the proof.
\end{proof}

\section{Slab theorems in arbitrary codimension}
Theorem \ref{thma} stated in the introduction describes the behavior of complete translating solitons contained in a slab under the geometric condition
$$
\beta_1+\beta_2\leq\frac{\pi}{2},
$$
where $\beta_1$ and $\beta_2$ are the principal angles between the tangent space and the fixed two-plane determined by the translating direction and the slab direction. We are now ready to prove it.

\begin{proofof}{Theorem \ref{thma}}
By assumption
$$
L\doteq\operatorname{osc}(\alpha)<\pi.
$$
After translating the immersion in the $W$-direction, we may assume that
$$
|\alpha|\leq\frac{L}{2}.
$$
Since $L<\pi$, we can choose $\sigma>1$ such that $\sigma L<\pi.$
Then
$$
|\sigma\alpha| \leq\frac{\sigma L}{2}<\frac{\pi}{2},
$$
and consequently
$$
\cos(\sigma\alpha)\geq \cos\left(\frac{\sigma L}{2}\right)>0.
$$

The proof will follow from the Omori-Yau maximum principle applied to the positive function $f_{\sigma}:M^m\to\R$ given by
$$
f_\sigma\doteq e^\y\cos(\sigma\alpha).
$$
Recall that from \eqref{thm:drift_identity} we have
$$\lap f_\sigma=\langle\nabla f_\sigma,\nabla\y
+\lambda \mathcal{X}\rangle-\lambda f_\sigma
+(1+\lambda)f_\sigma\big({\rm S}(\y,\alpha)-{\rm C}(\y,\alpha)\big),$$
where
$$
\mathcal{X}\doteq \nabla\y+\sigma\tan(\sigma\alpha)\nabla\alpha
\quad\text{and}\quad
\lambda=\frac{\sigma^2-1}{1+\sigma^2\tan^2(\sigma\alpha)}.
$$
The crucial quantity in \eqref{thm:drift_identity} is the
difference ${\rm S}(\y,\alpha)-{\rm C}(\y,\alpha)$.

Recall that
$$
{\rm C}(\y,\alpha)=|V^\top|^2|W^\top|^2-\langle V^\top,W^\top\rangle^2
=|V^\top\wedge W^\top|^2
$$
and
$$
{\rm S}(\y,\alpha)\doteq|V^\perp|^2|W^\perp|^2-\langle V^\perp,W^\perp\rangle^2=|V^\perp\wedge W^\perp|^2.
$$
The vectors $V^\top$ and $W^\top$ are the images of the orthonormal
basis $\{V,W\}$ of $P$ under the orthogonal projection $\pi_p^\top.$
Therefore $|V^\top\wedge W^\top|$
is the area distortion of $\pi_p^\top$. Since in terms of the principal angles $\beta_1$, $\beta_2$, its singular values are
given by $\cos\beta_1$ and $\cos\beta_2$, we obtain
\begin{equation}\label{eq:CS_principal_angles}
{\rm C}(\y,\alpha)=\cos^2\beta_1\cos^2\beta_2.
\end{equation}
Similarly
\begin{equation}\label{eq:N_principal_angles}
{\rm S}(\y,\alpha)= \sin^2\beta_1\sin^2\beta_2.
\end{equation}
By our assumption \eqref{b1b2intro} it follows that
$\cos(\beta_1+\beta_2)\geq0,$
and hence
$$
\cos\beta_1\cos\beta_2-\sin\beta_1\sin\beta_2\geq 0.
$$
Using \eqref{eq:CS_principal_angles} and
\eqref{eq:N_principal_angles}, we conclude that
\begin{equation}\label{eq:N_leq_CS}
{\rm S}(\y,\alpha)\leq{\rm C}(\y,\alpha).
\end{equation}
Therefore, since $f_\sigma>0$, we deduce from above that
\begin{equation}\label{eq:principal_angle_drift}
\lap f_\sigma
\leq
\langle\nabla f_\sigma,\nabla\y+\lambda \mathcal{X}\rangle
-\lambda f_\sigma.
\end{equation}
Now suppose, to the contrary, that $\y$ is bounded from below. Then there
exists $c_0\in\real{}$ such that $\y\geq c_0$ on $M^m$. Hence
$$
f_\sigma=e^\y\cos(\sigma\alpha)
\geq e^{c_0} \cos\left(\frac{\sigma L}{2}\right)>0.
$$
Thus
$$
\delta\doteq\inf_Mf_\sigma>0.
$$
By the properness of $M$ we may apply the Omori-Yau maximum principle to $f_\sigma$. Hence there exists a sequence
$\{p_j\}_{j\in\natural{}}\subset M$ such that
$$
f_\sigma(p_j)\to\delta,\quad
|\nabla f_\sigma(p_j)|\to 0\quad
\text{and}\quad
\liminf_{j\to\infty}\lap f_\sigma(p_j)\geq0.
$$
We next observe that all the coefficients in the first-order term of
\eqref{eq:principal_angle_drift} are uniformly bounded. Indeed,
by Lemma \ref{lem:grad_relations} we have
$$
|\nabla\y|=|V^\top|\leq1
\quad\text{and}\quad
|\nabla\alpha|=|W^\top|\leq1.
$$
Moreover,
$$
|\tan(\sigma\alpha)|
\leq
\tan\left(\frac{\sigma L}{2}\right).
$$
Consequently,
$$
|\mathcal{X}|\leq |\nabla\y|+\sigma|\tan(\sigma\alpha)|\,|\nabla\alpha|
\leq 1+ \sigma\tan\left(\frac{\sigma L}{2}\right).
$$
Thus $\mathcal X$ is uniformly bounded on $M$.
Furthermore,
$$
\lambda
=
\frac{\sigma^2-1}
{1+\sigma^2\tan^2(\sigma\alpha)}
$$
satisfies $0<\lambda\leq\sigma^2-1,$
while
$$
\lambda \geq \frac{\sigma^2-1}
{1+\sigma^2\tan^2\left(\frac{\sigma L}{2}\right)}\doteq\lambda_0>0.
$$
Evaluating \eqref{eq:principal_angle_drift} at $p_j$, we obtain
$$
\lap f_\sigma(p_j)
\leq\langle\nabla f_\sigma(p_j),\nabla\y(p_j)+\lambda(p_j)\mathcal{X}(p_j)
\rangle-\lambda_0f_\sigma(p_j),
$$
from where we arrive at the contradiction
$$
0\leq\liminf_{j\to\infty}\Delta f_\sigma(p_j)\leq\limsup_{j\to\infty}\Delta f_\sigma(p_j) \leq-\lambda_0\delta<0.
$$
Therefore $\y$ cannot be bounded from below and this concludes the proof.
\end{proofof}
We proceed with the proof of Theorem \ref{thmb}.
\begin{proofof}{Theorem \ref{thmb}}
It is well known that translating solitons can be viewed as minimal submanifolds with respect to a suitable analytic conformal change of the ambient metric; see, for example,
\cites{Ilmanen,MSS,Smoczyk12}. Consequently, translating solitons are real analytic.
By assumption
$$
\operatorname{osc}(\alpha)=\pi.
$$
After translating the immersion in the $W$-direction, we may assume without loss of generality that
$$
|\alpha|\leq\frac{\pi}{2},
$$
and consequently $\alpha_0=\sup\alpha-\pi/2=0$, and
$$
\cos(\alpha)\geq 0.
$$
Then by \eqref{thm:drift_identity} we have
\begin{equation}\label{f1}
\lap f=\langle\nabla f,\nabla\y\rangle
+f\big({\rm S}(\y,\alpha)-{\rm C}(\y,\alpha)\big),
\end{equation}
where by the principal angle assumption we conclude as in the proof of Theorem \ref{thma} that 
$${\rm S}(\y,\alpha)-{\rm C}(\y,\alpha)\le 0.$$
\textbf{(i)$\Rightarrow$(ii):} Suppose $f=e^\y\cos(\alpha)$ attains a local minimum. Since $f\ge 0$, all local minima must be non-negative. Since ${\rm S}(\y,\alpha)-{\rm C}(\y,\alpha)\le 0$ the strong elliptic maximum principle applied to \eqref{f1} shows that $f$ will be constant on any connected component of $M^m$ admitting a local non-negative minimum of $f$. Since we always assume that the translator is connected, we conclude that in such a case $f\equiv c$ for some constant $c\ge 0$ on all of $M^m$. Now $c>0$, because $c=0$ implies $\cos(\alpha)=0$ and then $\alpha\equiv\pm\pi/2$ which contradicts $\operatorname{osc}(\alpha)=\pi$. Therefore (i) implies (ii).

\textbf{(ii)$\Rightarrow$(iii):}  Suppose $f\equiv c>0$ for some constant $c$. Equation \eqref{f1} then implies that ${\rm S}(\y,\alpha)-{\rm C}(\y,\alpha)= 0$. 	Taking the logarithmic gradient of $f=c>0$, we obtain
$$
\nabla\y=\tan(\alpha)\nabla\alpha.
$$
Combining this with
$$
|\nabla\y|^2+|\nabla\alpha|^2
=
1+\mathrm C(\y,\alpha)-\mathrm S(\y,\alpha)
=
1,
$$
we obtain
$$
|\nabla\alpha|=\cos(\alpha),
$$
where we have used $\cos(\alpha)>0$. Define
$$
E\doteq\frac{\nabla\alpha}{\cos(\alpha)}.
$$
Then
$$
\nabla\alpha=\cos(\alpha)E
\quad\text{and}\quad
\nabla\y=\sin(\alpha)E.
$$
Equivalently, from Lemma \ref{lem:grad_relations},
$$
V^\top=\sin(\alpha)E
\quad\text{and}\quad
W^\top=\cos(\alpha)E.
$$

We claim that
$$
E=\sin(\alpha)V+\cos(\alpha)W.
$$
Indeed, the tangential component of the vector on the right-hand side is
$$
\sin(\alpha)V^\top+\cos(\alpha)W^\top=\bigl(\sin^2(\alpha)+\cos^2(\alpha)\bigr)E=E.
$$
Moreover,
$$
|\sin(\alpha)V+\cos(\alpha)W|=1=|E|,
$$
and hence its normal component vanishes.
Set
$$
N\doteq\cos(\alpha)V-\sin(\alpha)W.
$$
Then
$$
N^\top=\cos(\alpha)V^\top-\sin(\alpha)W^\top=0,
$$
so $N$ is a unit normal vector field. Since $V$ and $W$ are constant
ambient vector fields, for every $Y\in TM$ we have
$$
\overline\nabla_YE=Y(\alpha)\bigl(\cos(\alpha)V-\sin(\alpha)W\bigr)=
Y(\alpha)N.
$$
Thus $\overline\nabla_YE$ is normal. By the Gauss formula, we have
$$
\overline\nabla_YE=\nabla_YE+ A(Y,E),
$$
where $A$ is the second fundamental form. Therefore
$$
\nabla_YE=0
$$
for every $Y\in TM$. Since $E$ is a parallel unit vector field, the distributions
$\operatorname{span}\{E\}$ and $\mathcal D=E^{\perp}$
are parallel. Thus, by the de Rham decomposition theorem, $M$ is
locally a Riemannian product. In the present situation, the splitting
is in fact global.
Since $M$ is complete, every integral curve of $E$ is defined for all
time. Let $\gamma:\real{}\to M$ be an integral curve of $E$ and set
$$
a\doteq\alpha\circ\gamma.
$$
Because $\nabla\alpha=\cos(\alpha)E,$
we have
$$
a'=\cos a>0.
$$
Hence $a$ is strictly increasing and
$$
s-s_0=\int_{a(s_0)}^{a(s)}\frac{dt}{\cos t}.
$$
Since
$$
\int_{a(s_0)}^{\pi/2}\frac{dt}{\cos t}=\int_{-\pi/2}^{a(s_0)}\frac{dt}{\cos t}=+\infty,
$$
it follows that
$$
\lim_{s\to+\infty}a(s)=\frac{\pi}{2}\quad\text{and}\quad\lim_{s\to-\infty}a(s)=-\frac{\pi}{2}.
$$
Therefore
$$
(\alpha\circ\gamma)(\real{})=\left(-\frac{\pi}{2},\frac{\pi}{2}\right).
$$
In particular, every integral curve of $E$ intersects the level set
$$
\Sigma=\{\alpha=0\}
$$
at exactly one point. Since $|\nabla\alpha|=\cos\alpha=1$
on $\Sigma$, this level set is a smooth hypersurface of $M$.
Let $\{\Phi_s\}_{s\in\mathbb R}$ denote the flow of $E$. Since $E$ is parallel, it is a
Killing vector field and its flow preserves both the metric and the
distribution $\mathcal D$. Consequently, the map
$\Phi:\real{}\times\Sigma\to M$, given by
$$
\Phi(s,p)=\Phi_s(p),
$$
is a global isometry. Therefore,
$$
M\cong\real{}\times\Sigma
$$
as a Riemannian product. In particular, $\Sigma$ is complete.
We next show that this splitting is also extrinsic. Since
$$
\nabla\y=\sin(\alpha)E\quad\text{and}\quad\nabla\alpha=\cos(\alpha)E,
$$
for every $X\in\mathcal D$, we have
$$X(\y)=X(\alpha)=0.$$
Since $V$ and $W$ form an orthonormal basis of
$$P=\operatorname{span}\{V,W\},$$
the orthogonal projection of $F$ onto $P$ is
$$
\pi_P\circ F=\y V+\alpha W.
$$
It follows that the $P$-component of $F$ is constant along each leaf of
$\mathcal D$. On the other hand,
$$
dF(E)=E=\sin(\alpha)V+\cos(\alpha)W\in P,
$$
and hence
$$
d(\pi_{P^\perp}\circ F)(E)=0.
$$
Thus the $P^\perp$-component of $F$ is constant along the integral
curves of $E$.
Consequently, after an ambient translation, the immersion splits as
$$
F(s,p)=\Gamma(s)+F_\Sigma(p),
$$
where $\Gamma\subset P$ and $F_\Sigma(\Sigma)\subset P^\perp.$
In particular, identifying $P^\perp\simeq\real{m+k-2},$
we obtain an extrinsic product
$$
F(M)=\Gamma\times\Sigma.
$$
It remains to identify the curve $\Gamma$. Along an integral curve of
$E$, we have
$$
\frac{d\y}{ds}=E(\y)=\sin(\alpha)\quad\text{and}\quad\frac{d\alpha}{ds}
=E(\alpha)=\cos(\alpha).
$$
Therefore,
$$
\frac{d\y}{d\alpha}=\tan(\alpha),
$$
and hence
$$
\y=-\log\cos(\alpha)+\tilde c,
$$
for some constant $\tilde c$. Since $f=e^\y\cos(\alpha)\equiv c$, this implies $\tilde c=\log c$.
After translating in the $V$-direction, we may assume that $\tilde c=0$ and hence $c=1$.
Thus the $P$-component of the integral curves of $E$ is given by
$$
\Gamma(\alpha)=\alpha W-\log\cos(\alpha)\,V,
\quad\text{for}\,\,\,-\frac{\pi}{2}<\alpha<\frac{\pi}{2}.
$$
This is precisely the grim reaper curve in the plane $P$, translating
with unit speed in the direction $V$.
Finally, since the immersion is the orthogonal product $M=\Gamma\times\Sigma,$
its second fundamental form has no mixed terms, and its mean curvature
vector decomposes as
$$
\mathbf H=\mathbf H_\Gamma+\mathbf H_\Sigma.
$$
Since $V\in P$ and $T\Sigma\subset P^\perp$, the normal projection of
$V$ onto $M$ coincides with its normal projection onto $\Gamma$.
Moreover, since $\Gamma$ is the grim reaper translating in the
direction $V$,
$$
\mathbf H_\Gamma=V^\perp.
$$
Using the translator equation $\mathbf H=V^\perp,$
we conclude that $\mathbf H_\Sigma=0.$
Thus
$
\Sigma\subset P^\perp\simeq\real{m+k-2}
$
is a complete minimal submanifold.

\textbf{(iii)$\Rightarrow$(i):} Conversely, let us suppose that $M=\Gamma\times\Sigma,$
where $\Gamma\subset P$ is a grim reaper curve translating with unit
speed in the direction $V$ and $\Sigma\subset P^\perp$
is a complete minimal submanifold. After a translation in $P$, we may
write
$$
\Gamma(t)=tW-\log\cos(t)\,V,\quad\text{for}\,\,\,-\frac{\pi}{2}<t<\frac{\pi}{2}.
$$
Hence $\alpha=t$ and $\y=-\log\cos(t),$
and therefore
$$
f=e^\y\cos(\alpha)=1.
$$
Thus $f$ is constant on $M$ and, in particular, attains its positive infimum at every point.
\end{proofof}

\section{Applications}\label{applications}

In this section, we give several applications of
Theorem~\ref{thma}. As shown in
Section~\ref{Gauss}, the principal-angle condition \eqref{b1b2} is
automatically satisfied for hypersurfaces and coisotropic submanifolds,
while in the graphical setting it is equivalent to the weak
area-decreasing condition.

\subsection*{Codimension one.}\label{rem:codim_one_case}

As observed in Section~\ref{Gauss}, the principal-angle condition
\begin{equation}\label{b1b2}\tag{$\beta$}
\beta_1+\beta_2\leq\frac{\pi}{2},
\end{equation}
is automatically satisfied in codimension one.
Therefore, Theorem~\ref{thma} immediately gives the
following result.

\begin{corollary}
Let $F:M^m\to\real{m+1}$ be a properly immersed translating soliton moving in the unit direction $V$. Let
$W\in\real{m+1}$ be a constant unit vector orthogonal to $V$, and set
$\y\doteq\langle V,F\rangle$ and $\alpha\doteq\langle W,F\rangle.$
If $\operatorname{osc}(\alpha)<\pi$
then the height function $\y$ cannot be bounded from below.
\end{corollary}

\subsection*{Graphs}

We next consider graphical translators. Let $P=\operatorname{span}\{V,W\},$
and suppose that $M$ is given as the graph of a map $f:\Omega\subset P\to P^\perp.$
As shown in Section~\ref{Gauss}, the principal-angle condition
\eqref{b1b2} is equivalent to
$$
|\Lambda^2df|\leq1,
$$
that is, to $f$ being weakly area decreasing. Thus
Theorem~\ref{thma} applies naturally to this class
of graphical translators.

There are simple examples illustrating this situation. Let
$a\in\mathbb{R}\setminus\{0\}$ and consider the tilted grim reaper hyperplane in
$\mathbb{R}^3$
$$
F(x_1,x_2)=\Big(x_1,x_2,-(1+a^2)\log\Big(\cos\frac{x_1}{\sqrt{1+a^2}}\Big)+ax_2\Big),
$$
which translates with unit speed in the direction $V=e_3$; see
Subsection~\ref{tilted}. It is defined on
$$
\Big(-\frac{\pi}{2}\sqrt{1+a^2},\frac{\pi}{2}\sqrt{1+a^2}\Big)\times\mathbb{R}
$$
and is contained in a slab of width
$$
\pi\sqrt{1+a^2}>\pi
$$
in the direction $W=e_1$.
The same surface may be regarded as a graph over $P=\operatorname{span}\{e_1,e_3\}.$
Indeed,
$$
x_2=\frac{1}{a}\Big(x_3+(1+a^2)\log\Big(\cos\frac{x_1}{\sqrt{1+a^2}}\Big)\Big),
$$
and hence, setting
$$
\Omega=\Big(-\frac{\pi}{2}\sqrt{1+a^2},\frac{\pi}{2}\sqrt{1+a^2}\Big)\times\mathbb{R}\subset P,
$$
the tilted grim reaper hyperplane is the graph of $f:\Omega\to P^\perp=\operatorname{span}\{e_2\},$
given by
$$
f(x_1,x_3)=\frac{1}{a}\Big(x_3+(1+a^2)\log\Big(\cos\frac{x_1}{\sqrt{1+a^2}}\Big)\Big)e_2.
$$
Since $P^\perp$ is one-dimensional, $\operatorname{rank}df\leq1,$
and therefore $|\Lambda^2df|=0.$
Thus every nontrivially tilted grim reaper hyperplane in this family is a weakly
area-decreasing translating graph over $P$. Notice that these examples
are contained in slabs of width strictly greater than $\pi$.

The same examples may also be viewed as codimension-two translators in $\mathbb{R}^4$ by using the inclusion $\mathbb{R}^3\subset\mathbb{R}^4$ and taking
$P=\operatorname{span}\{e_1,e_3\}$ and $P^\perp=\operatorname{span}\{e_2,e_4\}.$
The corresponding graph has zero $e_4$-component and still satisfies
$\operatorname{rank}df\leq1.$ Hence it is weakly area decreasing. These examples are, however,
contained in the three-dimensional subspace
$$
\operatorname{span}\{e_1,e_2,e_3\}\subset\mathbb{R}^4,
$$
and therefore are not genuinely codimension-two examples.

In the graphical setting,
Theorem~\ref{thma} gives the following.

\begin{corollary}
Let $F:M^2\to\real{2+k}$ be a properly immersed translating soliton moving in the unit direction $V$. Let
$W\in\real{2+k}$ be a constant unit vector orthogonal to $V$, set $P=\operatorname{span}\{V,W\},$
and suppose that $M$ is the graph of a weakly area-decreasing map
$f:\Omega\subset P\to P^\perp.$
If
$$
\operatorname{osc}\langle W,F\rangle<\pi,
$$
then
$$
\inf_M\langle V,F\rangle=-\infty.
$$
In particular, there is no such complete translator if
$$
\Omega\subset\left\{x\in P:\langle V,x\rangle\geq c\,\,\,\text{and}\,\,\, a<\langle W,x\rangle<b\right\},
\quad b-a<\pi.
$$
\end{corollary}

\begin{proof}
By Section~\ref{Gauss}, the weak area-decreasing condition implies
\eqref{b1b2}. The first assertion therefore follows directly from
Theorem~\ref{thma}.
For the second assertion, since $f(x)\in P^\perp$ and $V,W\in P$, we have
$$
\langle V,F(x)\rangle=\langle V,x\rangle \quad\text{and}\quad \langle W,F(x)\rangle=\langle W,x\rangle.
$$
Thus, if
$$
\Omega\subset \left\{x\in P: \langle V,x\rangle\geq c\,\,\,\text{and}\,\,\, a<\langle W,x\rangle<b\right\},
$$
with $b-a<\pi$, then
$$
\langle V,F\rangle\geq c \quad\text{and}\quad \operatorname{osc}\langle W,F\rangle\leq b-a<\pi,
$$
contradicting the first assertion.
\end{proof}

\subsection*{Coisotropic and Lagrangian submanifolds.}

We finally consider coisotropic translators in $\complex{n}$. Recall
that an immersed submanifold $M^m\subset\complex{n}$ is coisotropic if
$$
J(N_pM)\subset T_pM,\quad\text{for any}\,\,\,p\in M,
$$
where $J$ denotes the standard complex structure.
Every Lagrangian submanifold is coisotropic; in the Lagrangian case, $m=n$ and
$$
J(N_pM)=T_pM, \quad\text{for any}\,\,\,p\in M.
$$
Let $F:M^m\to\complex{n}$ be a coisotropic translating soliton with
unit translating direction $V$ and set $W=JV.$
As proved in Section~\ref{Gauss}, the corresponding principal angles
satisfy
$$
\beta_1+\beta_2\leq{\pi}/{2}.
$$
In the Lagrangian case, equality holds:
$$
\beta_1+\beta_2={\pi}/{2}.
$$
Thus the geometric hypothesis of
Theorem~\ref{thma} is automatically satisfied in
the coisotropic setting. We obtain the following result, which extends
a theorem of Smoczyk \cite{Smoczyk19}*{Theorem~C} from the Lagrangian to
the coisotropic case.

\begin{corollary}
Let $F:M^m\to\complex{n}$ be a properly immersed coisotropic translating soliton moving in the unit direction $V$. Set $\y=\langle V,F\rangle$
and $\alpha=\langle JV,F\rangle.$  If $\operatorname{osc}(\alpha)<\pi$
then the translating height function $\y$ cannot be bounded from below.
\end{corollary}

Moreover, Theorem~\ref{thmb} extends a previous theorem of the third author \cite{Smoczyk19}*{Theorem~B}, proved in the Lagrangian setting.

\section{Ancient solutions confined in slabs}
Throughout this section we will always assume that one can apply the parabolic Omori-Yau maximum principle to the mean curvature flow under consideration.	To prove Theorem~\ref{thmc} we need several ingredients. We start with the computation of a crucial evolution equation.
\begin{lemma}\label{evol fsigma}
	Let $F: M^m \times [t_0, T) \to \R^{m+k}$ be evolving under mean curvature flow $\dt F = \mathbf H$, where $t_0 < T$ is arbitrary. Let $V, W \in \R^{m+k}$ be constant orthonormal vectors, and let $u = \la V, F \ra$ and $\alpha = \la W, F \ra$ be the corresponding coordinate functions. For any $\sigma \in \R$, define $f_\sigma = e^{u-t} \cos(\sigma \alpha)$. Then $f_\sigma$ satisfies the evolution equation
	\begin{equation}\label{eq:f_sigma_evolution}
		(\dt - \De + 2\la \na u, \na \cdot \ra) f_\sigma = f_\sigma \left( -1 + |\na u|^2 + \sigma^2 |\na \alpha|^2 \right).
	\end{equation}
	Moreover, at any point where $\tan(\sigma\alpha)$ is well defined, we have
		\begin{equation}\label{eq:f_sigma_evolution_2}
		\left(\partial_t - \Delta + \langle Y, \nabla \cdot \rangle\right) f_\sigma = \lambda f_\sigma + (1 + \lambda) f_\sigma \left( \mathrm{C}(u, \alpha) - \mathrm{S}(u, \alpha) \right),
	\end{equation}
	where $\mathrm C(u, \alpha)$ and $\mathrm S(u, \alpha)$ are defined as in Lemma \ref{lem:grad_relations}, and where
    $$\lambda \doteq \frac{\sigma^2 - 1}{1 + \sigma^2 \tan^2(\sigma \alpha)}$$
    and $Y \doteq 2\nabla u + \lambda \mathcal{X}$ with $\mathcal{X} \doteq \nabla u + \sigma \tan(\sigma \alpha) \nabla \alpha$.
\end{lemma}
\begin{proof}
Under mean curvature flow
$$\dt F = \mathbf H=\Delta F,$$
the coordinate functions satisfy
$$(\dt - \De) u = 0\quad\text{and}\quad (\dt - \De) \alpha = 0.$$
Computing the spatial gradient of $f_\sigma = e^{u-t} \cos(\sigma \alpha)$ gives
	\begin{equation}\label{eq:grad_f_neg}
		\na f_\sigma = f_\sigma \na u - \sigma e^{u-t} \sin(\sigma \alpha) \na \alpha.
	\end{equation}
Taking the inner product of \eqref{eq:grad_f_neg} with $\na u$ implies
	\[
	\la \na f_\sigma, \na u \ra = f_\sigma |\na u|^2 - \sigma e^{u-t} \sin(\sigma \alpha) \la \na u, \na \alpha \ra,
	\]
which allows us to express the mixed term as
	\begin{equation}\label{eq:mixed_term_f_neg}
		\sigma e^{u-t} \sin(\sigma \alpha) \la \na u, \na \alpha \ra = f_\sigma |\na u|^2 - \la \na f_\sigma, \na u \ra.
	\end{equation}
Next, we calculate the time derivative of $f_\sigma$.
\begin{eqnarray*}
		\dt f_\sigma &=& -e^{u-t} \cos(\sigma \alpha) + e^{u-t} \cos(\sigma \alpha) \dt u - \sigma e^{u-t} \sin(\sigma \alpha) \dt \alpha \\
		&=& -f_\sigma + f_\sigma \De u - \sigma e^{u-t} \sin(\sigma \alpha) \De \alpha.
\end{eqnarray*}
Applying the Laplacian operator $\De$ to $f_\sigma$ gives
	\begin{eqnarray*}
		\De f_\sigma &=& f_\sigma \De u - \sigma e^{u-t} \sin(\sigma \alpha) \De \alpha + f_\sigma |\na u|^2 - \sigma^2 f_\sigma |\na \alpha|^2 \\
		&\quad& - 2\sigma e^{u-t} \sin(\sigma \alpha) \la \na u, \na \alpha \ra.
	\end{eqnarray*}
Then, subtracting $\De f_\sigma$ from $\dt f_\sigma$, the terms involving $\De u$ and $\De \alpha$ cancel, and we obtain
	\[
	(\dt - \De) f_\sigma = -f_\sigma - f_\sigma |\na u|^2 + \sigma^2 f_\sigma |\na \alpha|^2 + 2\sigma e^{u-t} \sin(\sigma \alpha) \la \na u, \na \alpha \ra.
	\]
Substituting \eqref{eq:mixed_term_f_neg} into the last term yields
	\begin{eqnarray*}
		(\dt - \De) f_\sigma &=& -f_\sigma - f_\sigma |\na u|^2 + \sigma^2 f_\sigma |\na \alpha|^2 + 2\left( f_\sigma |\na u|^2 - \la \na f_\sigma, \na u \ra \right) \\
		&=& -f_\sigma + f_\sigma |\na u|^2 + \sigma^2 f_\sigma |\na \alpha|^2 - 2\la \na f_\sigma, \na u \ra.
	\end{eqnarray*}
Rearranging terms to isolate the drift operator gives \eqref{eq:f_sigma_evolution}. 
Defining
$$\lambda \doteq \frac{\sigma^2 - 1}{1 + \sigma^2 \tan^2(\sigma \alpha)},$$
we observe
$$1 + \lambda = \frac{\sigma^2 (1 + \tan^2(\sigma \alpha))}{1 + \sigma^2 \tan^2(\sigma \alpha)}.$$ 
Since
\[
-1 + |\nabla u|^2 + \sigma^2 |\nabla \alpha|^2 = (1 + \lambda)(-1 + |\nabla u|^2 + |\nabla \alpha|^2) + \lambda - \lambda \left(|\nabla u|^2 - \sigma^2 \tan^2(\sigma \alpha) |\nabla \alpha|^2\right),
\]
and taking into account that
$$-1+ |\nabla u|^2 + |\nabla \alpha|^2=\mathrm{C}(u, \alpha)-\mathrm{S}(u, \alpha)$$
we can rewrite equation \eqref{eq:f_sigma_evolution} in the form
\begin{eqnarray*}
		\left(\partial_t - \Delta + 2\langle \na u, \nabla \cdot \rangle\right) f_\sigma &=& \lambda f_\sigma + (1 + \lambda) f_\sigma \left( \mathrm{C}(u, \alpha)- \mathrm{S}(u, \alpha) \right)\\
        &&-\lambda f_\sigma\left(|\nabla u|^2 - \sigma^2 \tan^2(\sigma \alpha) |\nabla \alpha|^2\right).
\end{eqnarray*}
Noting that
$$\nabla \log f_\sigma = \nabla u - \sigma \tan(\sigma \alpha) \nabla \alpha,$$
taking the inner product with
$$\mathcal{X} = \nabla u + \sigma \tan(\sigma \alpha) \nabla \alpha$$
gives
\[
\langle \nabla \log f_\sigma, \mathcal{X} \rangle = |\nabla u|^2 - \sigma^2 \tan^2(\sigma \alpha) |\nabla \alpha|^2,
\]
which implies
$$
-\lambda f_\sigma \left(|\nabla u|^2 - \sigma^2 \tan^2(\sigma \alpha) |\nabla \alpha|^2\right)= -\langle \nabla f_\sigma, \lambda \mathcal{X} \rangle.
$$
Substituting these identities back into the evolution equation gives:
\[
(\partial_t - \Delta + 2\langle \nabla u, \nabla \cdot \rangle +\lambda \langle \mathcal{X}, \nabla \cdot \rangle) f_\sigma = \lambda f_\sigma + (1 + \lambda) f_\sigma \left( \mathrm{C}(u, \alpha) - \mathrm{S}(u, \alpha) \right).
\]
Combining the drift terms into a single vector field
$$Y \doteq 2\nabla u + \lambda \mathcal{X},$$
we obtain \eqref{eq:f_sigma_evolution_2}.
\end{proof}
\begin{lemma}\label{lemma unbounded}
	Let $F: M^m \times [t_0, T) \to \mathbb{R}^{m+k}$ be a family of immersions evolving under the mean curvature flow $\partial_t F = \mathbf{H}$.
   Assume that, for every $s\in(t_0,T)$, the spacetime immersion
   $(x,t)\mapsto(F(x,t),t)$ is proper on $M\times[t_0,s]$.
    Let $V, W \in \mathbb{R}^{m+k}$ be constant orthonormal vectors, and set $u \coloneqq \langle V, F \rangle$ and $\alpha \coloneqq \langle W, F \rangle$. Assume that the principal angles between the plane $P \coloneqq \operatorname{span}\{V, W\}$ and the tangent spaces of the solution satisfy the angle condition \eqref{b1b2intro} on $M \times [t_0, T)$. 
    Suppose that  $\operatorname{osc}(\alpha) < \pi$ so that, without loss of generality, there exist $\sigma > 1$ and $\varepsilon_0 > 0$ such that $\cos(\sigma\alpha) \ge \varepsilon_0$. Then $f_\sigma = e^{u-t} \cos(\sigma \alpha)$ cannot attain a positive infimum in $M^m \times (t_0, T)$. That is, if 
	\[
	I \doteq \inf_{M^m \times [t_0, T)} f_\sigma > 0,
	\]
	then $I = \inf_{M^m} f_\sigma(\cdot, t_0)$.
\end{lemma}
\begin{proof}
	Suppose for the sake of contradiction that 
	\[
	0<I = \inf_{M^m \times [t_0, T)} f_\sigma < \inf_{M^m} f_\sigma(\cdot, t_0).
	\]
	Since $\sigma > 1$ and $\cos(\sigma \alpha) \ge \varepsilon_0 > 0$, the function
    $$\lambda(x,t) = \frac{\sigma^2 - 1}{1 + \sigma^2 \tan^2(\sigma \alpha)}$$
    is uniformly bounded and strictly positive, i.e.,
    $$0 < \lambda_0 \le \lambda(x, t) \le \sigma^2 - 1.$$
    Furthermore, since $\mathrm{C}(u, \alpha) - \mathrm{S}(u, \alpha) \ge 0$ and $1 + \lambda > 0$, we have:
	\[
	(1 + \lambda) f_\sigma \left( \mathrm{C}(u, \alpha) - \mathrm{S}(u, \alpha) \right) \ge 0.
	\]
	Thus, $f_\sigma$ satisfies the differential inequality:
	\[
	\partial_t f_\sigma - \Delta f_\sigma + \langle Y, \nabla f_\sigma \rangle \ge \lambda f_\sigma \ge \lambda_0 f_\sigma.
	\]
	Because $I < \inf_{M^m} f_\sigma(\cdot, t_0)$, the infimum $I$ cannot be realized or approached on the initial boundary $M^m \times \{t_0\}$. Thus, by the generalized parabolic Omori-Yau maximum principle, there exists a sequence $(x_j, t_j) \in M^m \times (t_0, T)$ with $t_j > t_0$ such that:
    \begin{enumerate}[(1)]
		\item $\lim f_\sigma(x_j, t_j) = I$,
        \medskip
		\item $\lim |\nabla f_\sigma(x_j, t_j)| = 0$,
        \medskip
		\item $\limsup (\partial_t-\Delta)f_\sigma(x_j, t_j) \le 0$.
	\end{enumerate}
    We also observe that the vector field $Y$ is uniformly bounded. Indeed, $|\nabla u|\leq1$ and $|\nabla\alpha|\leq1$, while the assumption
    $\cos(\sigma\alpha)\geq\varepsilon_0>0$ implies that $|\tan(\sigma\alpha)|$ is uniformly bounded. Since $0<\lambda\leq\sigma^2-1$, it follows that
    $Y=2\nabla u+\lambda X$ is uniformly bounded.
	Evaluating the differential inequality along this sequence yields:
	\begin{eqnarray*}
		0 &\ge& \limsup_{j \to \infty} \left( \partial_t f_\sigma - \Delta f_\sigma + \langle Y, \nabla f_\sigma \rangle \right)(x_j, t_j) \\
		&\ge& \limsup_{j \to \infty} (\lambda_0 f_\sigma(x_j, t_j))\\
        &=& \lambda_0 I.
	\end{eqnarray*}
	Since $\lambda_0 > 0$ and $I > 0$, we have $\lambda_0 I > 0$, which implies $0 \ge \lambda_0 I > 0$, a contradiction. Hence, our assumption that $I < \inf_{M^m} f_\sigma(\cdot, t_0)$ must be false, establishing that
    $$I = \inf_{M^m} f_\sigma(\cdot, t_0).$$
    This completes the proof.
\end{proof}
\begin{proofof}{Theorem \ref{thmc}}
Since $F:M^m\times(-\infty,T)\to\real{m+k}$ is an ancient solution, we may assume without loss of generality that $T<\infty$. By assumption, we have $\operatorname{osc}(\alpha) < L$ for some width $L > 0$. By translating the solution $F$ if necessary in the direction of $W$, we may assume that 
$$\vert\inf\alpha\vert = \sup\alpha < \frac{L}{2}.$$ 

We divide the proof into several steps. 
\begin{enumerate}[(a)]
    \item Let us first assume that $L \le \pi$. Then, since $\vert\inf\alpha\vert = \sup\alpha < \pi/2$, we can find constants $\sigma > 1$, $\varepsilon_0 > 0$ such that $\cos(\sigma\alpha) \ge \varepsilon_0 > 0$ on $M^m \times (-\infty, T)$. We prove the unboundedness of $u$ from below by contradiction. Assume $u(x, t) \ge C_0$ for all $(x, t) \in M^m \times (-\infty, T)$. For any $t_0 \in (-\infty, T)$, consider
    $$f_\sigma(x, t) = e^{u(x, t)-t} \cos(\sigma \alpha(x, t)).$$
    In view of $u \ge C_0$ and $\cos(\sigma \alpha) \ge \varepsilon_0$, we have
    $$f_\sigma(x, t) \ge \varepsilon_0 e^{C_0} e^{-t} > 0.$$
	Applying Lemma \ref{lemma unbounded} on $M^m \times [t_0, T)$, since
    $$\inf_{M^m \times [t_0, T)} f_\sigma \ge \varepsilon_0 e^{C_0} e^{-T} > 0,$$
    the spatial-temporal infimum must be realized at $t = t_0$:
	\[
	\inf_{M^m \times [t_0, T)} f_\sigma = \inf_{M^m} f_\sigma(\cdot, t_0) = e^{-t_0} \inf_{M^m} \big( e^u \cos(\sigma \alpha) \big)(\cdot, t_0).
	\]
	Fixing a target time $t_1 \in (t_0, T)$, note that $M^m \times [t_1, T) \subset M^m \times [t_0, T)$, so
	\begin{eqnarray*}
	e^{-t_1} \inf_{M^m} \big( e^u \cos(\sigma \alpha) \big)(\cdot, t_1) 
	&=& \inf_{M^m} f_\sigma(\cdot, t_1) \\
	&\ge& \inf_{M^m \times [t_0, T)} f_\sigma = e^{-t_0} \inf_{M^m} \big( e^u \cos(\sigma \alpha) \big)(\cdot, t_0).
	\end{eqnarray*}
	Since $u \ge C_0$ and $\cos(\sigma \alpha) \ge \varepsilon_0$, we have
    $$\inf_{M^m} \big( e^u \cos(\sigma \alpha) \big)(\cdot, t_0) \ge \varepsilon_0 e^{C_0},$$ which yields
	\[
	\inf_{M^m} \big( e^u \cos(\sigma \alpha) \big)(\cdot, t_1) \ge \varepsilon_0 e^{C_0} e^{t_1 - t_0}.
	\]
	Taking $t_0 \to -\infty$, the factor $e^{t_1 - t_0} \to +\infty$, which forces $$\inf_{M^m} \big( e^u \cos(\sigma \alpha) \big)(\cdot, t_1) = +\infty,$$
	a contradiction. Thus, $u$ cannot be uniformly bounded from below.
    
    \item Set $\tilde{V} = -V$, so
    $\tilde{u} = \la \tilde{V}, F \ra = -u$. Since
    $$\mathrm C(\tilde{u}, \alpha) - \mathrm S(\tilde{u}, \alpha) = \mathrm C(u, \alpha) - \mathrm S(u, \alpha) \ge 0,$$
    applying part (a) to $\tilde{u}$ gives $\inf \tilde{u} = -\infty$, which implies $\sup u = +\infty$. Therefore $u$ can be neither uniformly bounded from above nor from below on $M^m\times(-\infty,T)$.
    \item Now suppose $L>\pi$ and choose $\lambda>L/\pi$. Define the homothetically rescaled solution $$\tilde{F}(x, s) \doteq \frac{1}{\lambda} F(x, \lambda^2 s)$$ 
	on $M^m \times (-\infty, T/\lambda^2)$. Since the slab is centered at 0, the rescaled coordinate function $\tilde{\alpha} = \la W, \tilde{F} \ra$ takes values in $[-{L}/{2\lambda}, {L}/{2\lambda}]$ of width $\tilde{L} = L/\lambda < \pi$. Applying parts (a) and (b) to $\tilde{F}$ establishes that $\tilde{u}_V$ is completely unbounded, whence $u_V = \lambda \tilde{u}_V$ is also completely unbounded.
    \item Now let $F: M^m \times (-\infty, T) \to \mathbb{R}^{m+1}$ be an ancient hypersurface evolving under the mean curvature flow. 
    Then, the principal angle condition holds for any choice of $V$ orthogonal to $W$. 
    Consequently, it follows immediately from the above that all coordinate functions $u_V$ are unbounded from both above and below on $M^m \times (-\infty, T)$.
\end{enumerate}
This completes the proof. 
\end{proofof}
As a corollary we obtain the following statement, which also follows from \cite{CMancient}*{Theorem 5}.
\begin{corollary}[Non-existence of Ancient Solutions Confined to Cylinders]
	In codimension-one, a direct geometric consequence of Theorem~\ref{thmc} is the non-existence\footnote{Note that we assume the applicability of the Omori-Yau maximum principle throughout this section. There exist complete immersed minimal surfaces contained in a ball and hence also in a cylinder, as shown by Nadirashvili~\cite{nadi}. This does not contradict Theorem~\ref{thmc} because these examples do not allow the application of the Omori-Yau maximum principle; in particular, they are not properly immersed.} of ancient solutions confined within a solid cylinder 
\[
\mathcal{C}_R = \{ x \in \mathbb{R}^{m+1} : \operatorname{dist}(x, L_1) \le R \}
\]
around a 1-dimensional axis $L_1$ for $m \ge 2$. 
Indeed, no such complete hypersurface (whether static or evolving) can exist inside the cylinder $\mathcal{C}_R$ when $m \ge 2$, because Theorem~\ref{thmc} implies unboundedness along all directions orthogonal to any slab direction, whereas the cylinder bounds the hypersurface in $m$ independent directions.
In contrast, for $m = 1$ (curves in $\mathbb{R}^2$), the cylinder $\mathcal{C}_R$ reduces to an infinite strip, and besides static straight lines parallel to $L_1$ there do exist complete ancient solutions contained in $\mathcal{C}_R$, e.g. the scaled grim reaper and the paper-clip solution (see Figure \ref{fig: AngenentOval}).

\end{corollary}

\begin{remark}[The Ancient Pancake Solution]
	A central explicit model that illustrates Theorem~\ref{thmc} is the paper-clip solution, see Figure \ref{fig: AngenentOval},
    together with its rotationally symmetric higher-dimensional generalization, the \emph{ancient pancake solution}.
	The paper-clip is the compact ancient curve in $\mathbb{R}^2$ confined to a slab, given by the implicit equation
	\[
	\cos(\alpha) = e^t \cosh(u), \quad t \in (-\infty, 0).
	\]
    For $m > 1$, the ancient pancake solution $M_t^m \subset \mathbb{R}^{m+1}$ is a compact, convex, symmetric $O(1) \times O(m)$ solution that is also contained in a slab for $t \to -\infty$. Although it lacks an elementary closed-form formula, it exhibits the same qualitative behavior: as $t \to -\infty$, the width of the slab remains bounded while all transverse height functions $u_V$ expand to cover $(-\infty, +\infty)$, thus realizing the behavior in Theorem~\ref{thmc}.
	\image{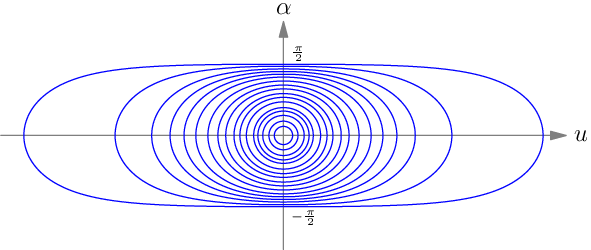}{The paper-clip solution.}{fig: AngenentOval}{.6}
\end{remark}

\bibliographystyle{y2k}

\end{document}